\documentclass [11pt]{amsart}
\usepackage {amsmath, amsthm, amssymb, amscd, mathrsfs, graphicx, stmaryrd, extarrows, color, url, enumerate, epsf, hyperref}
\usepackage[all,knot,arc]{xy}
\usepackage[text={6.5in,9in},centering,letterpaper,dvips]{geometry}

\newtheorem {theorem}{Theorem}[section]

\newtheorem {corollary}[theorem]{Corollary}

\theoremstyle{definition}
\newtheorem {definition}[theorem]{Definition}
\theoremstyle{remark}
\newtheorem {remark}[theorem]{Remark}
\newtheorem {example}[theorem]{Example}
\newtheorem {fact}[theorem]{Fact}

\newcommand\eps{\varepsilon}

\newcommand\Z{\mathbb{Z}}
\newcommand\R{\mathbb{R}}
\newcommand\Q{\mathbb{Q}}

\newcommand\T{\mathbb{T}}
\newcommand\Ta{\mathbb{T}_\alpha}
\newcommand\Tb{\mathbb{T}_\beta}

\newcommand\Xs{\mathbb{X}}
\newcommand\Os{\mathbb{O}}

\newcommand\Rect{\mathrm{Rect}}

\def\Sym{\mathrm{Sym}}

\def\Os{\mathbb O}

\def\zz {{\mathbb{Z}}}
\def\rr {{\mathbb{R}}}

\def\del {{\partial}}
\def\tt {{\mathfrak{t}}}
\def\ss {{\mathfrak{s}}}

\def\spinc {{\operatorname{Spin^c}}}

\def\kr{\operatorname{KR}}
\def\rk {{\operatorname{rank}}}

\def\fin\qedhere
\def\pr {{\text{pr}}}

\def\from {{\leftarrow}}

\def\M {\mathcal {M}}
\def\Mh {\widehat{\M}}

\def\x{\mathbf x}
\def\y{\mathbf y}

\def\Ring {\mathcal R}
\def\K {\mathcal K}

\newcommand\alphas{\boldsymbol\alpha}
\newcommand\betas{\boldsymbol\beta}

\def\ws{\mathbf w}
\def\zs{\mathbf z}

\newcommand\EmptyRect{\Rect^\circ}

\def\F {{\mathcal{F}}}

\def\He {\mathcal{H}}

\def\fin\qedhere
\def\from {{\leftarrow}}
\def\CF {\mathit{CF}}
\def\HF {\mathit{HF}}

\newcommand\HFplus {\HF^+}
\newcommand \CFplus {\CF^+}

\newcommand \CFminus {\CF^-}
\newcommand \HFminus {\HF^-}

\newcommand \HFinf {\HF^{\infty}}

\def\A{\mathcal{A}}
\def\Tor{\operatorname{Tor}}
\def\sl {\mathfrak{sl}}
\def\gl {\mathfrak{gl}}
\def\CFinfty {\CF^\infty}
\def\HFinfty {\HF^\infty}
\def\T{\mathcal{T}}
\def\XX{\mathbb{X}}

\newcommand\HFhatold{\widehat{\HF}}

\newcommand\HFhat{\widehat{\mathit{HF}}}
\newcommand\CFhat{\widehat{\mathit{CF}}}

\newcommand\HFLm{\mathit{HFL}^-}
\newcommand\CFLm{\mathit{CFL}^-}

\def\gCFKt{\widetilde{\mathit{gCFK}}}
\def\gCFKhat{\widehat{\mathit{gCFK}}}
\def\gCFKm{\mathit{gCFK}^-}
\def\HFMinus{\mathbf{HF}^-}

\def\In{\operatorname{In}}
\def\Out{\operatorname{Out}}

\def\HFLhat{\widehat{{HFL}}}
\def\CFKhat{\widehat{{CFK}}}
\def\CFKt{\widetilde{CFK}}
\def\CFKp{\mathit{CFK}^+}
\def\HFKt{\widetilde{HFK}}
\def\CFKm{\mathit{CFK}^-}
\def\CFKi{\mathit{CFK}^{\infty}}
\def\HFKhat{\widehat{\mathit{HFK}}}
\def\HFKm{\mathit{HFK}^-}
\def\CFKinfty{\CFKi}

\begin{document}

\title{An introduction to knot Floer homology}
\author[Ciprian Manolescu]{Ciprian Manolescu}
\thanks {The author was partially supported by NSF grant number DMS-1104406.}
\address {Department of Mathematics, UCLA, 520 Portola Plaza\\ 
Los Angeles, CA 90095}
\email {cm@math.ucla.edu}
\begin {abstract}
This is a survey article about knot Floer homology. We present three constructions of this invariant: the original one using holomorphic disks, a combinatorial description using grid diagrams, and a combinatorial description in terms of the cube of resolutions. We discuss the geometric information carried by knot Floer homology, and the connection to three- and four-dimensional topology via surgery formulas. We also describe some conjectural relations to Khovanov-Rozansky homology.
\end {abstract}

\maketitle

\section {Introduction}
Knot Floer homology is an invariant of knots and links in three-manifolds. It was introduced independently by Ozsv\'ath-Szab\'o \cite{Knots} and Rasmussen \cite{RasmussenThesis} around 2002. Since then it has grown into a large subject. Its importance lies in the fact that it contains information about several non-trivial geometric properties of the knot (genus, slice genus, fiberedness, effects of surgery, etc.) Furthermore, knot Floer homology is computable: There exist general algorithms that can calculate it for arbitrary knots. These algorithms tend to get slow as the complexity of the knot increases, but there are also different methods that can be applied to special classes of knots and give explicit answers.

The origins of knot Floer homology lie in gauge theory or, more precisely, in the Seiberg-Witten equations \cite{SW1, SW2, Witten}. These equations play an important role in four- and three-dimensional topology. In particular, given a  three-dimensional manifold $Y$, by studying the equations on $\R \times Y$, one can construct an invariant called the Seiberg-Witten Floer homology of $Y$ \cite{MarcolliWang, Spectrum, FroyshovSW, KMbook}. Inspired by the Atiyah-Floer conjecture \cite{AtiyahFloer}, Ozsv\'ath and Szab\'o developed Heegaard Floer theory as a symplectic geometric  replacement for gauge theory. In \cite{HolDisk}, \cite{HolDiskTwo}, they used Gromov's theory of pseudo-holomorphic curves to construct an invariant of closed 3-manifolds called Heegaard Floer homology. Knot Floer homology is a relative version of Heegaard Floer homology, associated to a pair consisting of a $3$-manifold and a null-homologous knot in it. It is now known that Heegaard Floer homology and Seiberg-Witten Floer homology are isomorphic \cite{CGH2, KLT1}. Thus, knot Floer homology can be thought of as encoding something about the Seiberg-Witten equations on $\R$ times the knot complement. This is only a heuristic, since at the moment no such direct gauge theoretic description exists. Nevertheless, some Seiberg-Witten counterparts to knot Floer homology were constructed in \cite{KMsutures, Kutluhan}.

Knot Floer homology is very similar in structure to knot homologies coming from representation theory, such as those introduced by Khovanov \cite{Khovanov} and Khovanov-Rozansky \cite{KR1, KR2}. There are three (mostly conjectural) relations between knot Floer homology and the Khovanov-Rozansky theories:
\begin{itemize} 
\item
Just as the Khovanov-Rozansky homologies are associated to the standard representation of the quantum group $U_q\bigl( \sl(n) \bigr)$, knot Floer homology is supposed to correspond to the standard representation of $U_q \bigl( \gl(1 | 1) \bigr)$. This connection has not yet been made explicit (except at the level of Euler characteristics), but there are various hints in this direction \cite{DMcornered, Tian1};
\item
There are structural parallels to Khovanov's $\sl(2)$ homology, discussed in Rasmussen's survey \cite{RasSurvey}. In particular, in all observed examples, Khovanov homology has rank at least as large as that of knot Floer homology, which may indicate that there is a spectral sequence connecting the two theories;
\item Dunfield, Gukov and Rasmussen \cite{DGR} have conjectured the existence of a $d_0$ differential (or, more realistically, a spectral sequence) starting at the Khovanov-Rozansky triply graded categorification of the HOMFLY-PT polynomial, and converging to knot Floer homology. 
\end{itemize}

The purpose of this survey is to give a general introduction to knot Floer homology, and to serve as a (necessarily incomplete) guide to the literature. The article is aimed at graduate students and researchers, especially those in related fields. 

Rather than first giving the definition of knot Floer homology, we will start by advertising some of its major  properties and applications; this is done in Section~\ref{sec:appl}. For example, we will discuss to what extent knot Floer homology distinguishes knots from each other.

In Section~\ref{sec:knot} we sketch the original definition of knot Floer homology, following \cite{Knots, RasmussenThesis}. Just like the Heegaard Floer homology of three-manifolds, knot Floer homology was first defined by counting pseudo-holomorphic curves in the symmetric product of a Heegaard surface. The appearance of the symmetric product is natural in view of the gauge theoretic origins of the subject. Points on the symmetric product describe solutions to the vortex equations on the surface, and the vortex equations are the two-dimensional reduction of the Seiberg-Witten equations.

This first  definition of knot Floer homology is the most flexible and the most useful for establishing the various properties of the invariant. Its main drawbacks are that it requires some familiarity with symplectic geometry, and that it is not combinatorial. (There is no known general algorithm for counting pseudo-holomorphic curves.) By making suitable choices of Heegaard diagrams (and thus making the pseudo-holomorphic curve counts more tractable) and/or by making use of the various properties of knot Floer homology, one can give several alternative, fully combinatorial definitions---at least for the case of knots in $S^3$. By now, there are several such constructions in the literature:
\begin{enumerate}
\item one by Manolescu-Ozsv\'ath-Sarkar, using grid diagrams \cite{MOS}; cf. also \cite{MOST};
\item another by Sarkar and Wang, using nice diagrams \cite{SarkarWang};  
\item another by Ozsv\'ath and Szab\'o, using a cube of resolutions \cite{CubeResolutions};
\item another by Baldwin and Levine, in terms of spanning trees \cite{BaldwinLevine};
\item yet another recently announced by Ozsv\'ath and Szab\'o, based on bordered Floer homology.
\end{enumerate}

In this paper we will only discuss (1) and (3). The grid diagram construction is conceptually the simplest, and is described in Section~\ref{sec:grids}. The cube of resolutions construction is the one closest in spirit to the Khovanov-Rozansky homologies, and can be used as a starting point for exploring the Dunfield-Gukov-Rasmussen conjecture; we will discuss it in Section~\ref{sec:cube}. 

Lastly, in Section~\ref{sec:surgery} we will outline how knot Floer homology fits into the Heegaard Floer theory developed by Ozsv\'ath and Szab\'o. The Heegaard Floer homology of a three-manifold obtained by surgery on a knot can be computed in terms of the knot Floer complex, using surgery formulas \cite{IntSurg, RatSurg}. These formulas can be extended to surgeries on links, and to link presentations of four-manifolds \cite{LinkSurg}. One may hope that the Heegaard Floer surgery formulas could serve as a model for extending Khovanov-Rozansky homology to three- and four-manifolds; we will explain the difficulties inherent in such a program.

\medskip
 \textbf {Acknowledgments.} This survey is based on lectures given by the author at the Summer School on the Physics and Mathematics of Link Homology, held at Montr{\'e}al in June-July 2013. The author would like to thank the participants and the organizers for making this event possible, and for their interest in the subject. Thanks go also to Jennifer Hom, \c{C}a\u{g}atay Kutluhan and Yajing Liu for helpful comments on an earlier version of the article.
 
 \section {Properties and applications}
\label {sec:appl}

Knot Floer homology can be defined for null-homologous knots in arbitrary three-manifolds. However, for simplicity, in this paper we focus on knots in the three-sphere.

\subsection{General form} \label{sec:general} Let $K \subset S^3$ be an oriented knot. There are several different variants of the knot Floer homology of $K$. The simplest is the hat version, which takes the form of a bi-graded, finitely generated Abelian group
$$\HFKhat(K) = \bigoplus_{i, s \in \Z} \HFKhat_i(K, s).$$
Here, $i$ is called the {\em Maslov} (or {\em homological}) grading, and $s$ is called the Alexander grading. The graded Euler characteristic of $\HFKhat$ is the Alexander-Conway polynomial:
\begin{equation}
\label {eq:alex}
 \sum_{s, i\in \zz} (-1)^i q^s \cdot \rk_{\Z}  \bigl( \HFKhat_i(K, s) \bigr) = \Delta_K(q).
 \end {equation}

Another version of knot Floer homology, called {\em minus} and denoted $\HFKm$, has the form of a bi-graded module over the polynomial ring $\Z[U]$, and contains more information than $\HFKhat$.  The most complete version, which has even more information, is not really a homology group but rather a doubly-filtered chain complex denoted $\CFKi$, well-defined up to filtered chain homotopy equivalence. We call $\CFKi$ the {\em full knot Floer complex}. These and other related variants are discussed in Section~\ref{sec:knot}. 

Knot Floer homology can be extended to links. In that setting, the Alexander and Maslov gradings may take half-integer values, and the Euler characteristic of the theory is the Alexander-Conway polynomial multiplied by the factor $(q^{-1/2} - q^{1/2})^{\ell -1}$, where $\ell$ is the number of components of the link. Furthermore, there is a refinement of knot Floer homology called {\em link Floer homology} \cite{Links}, which admits an Alexander multi-grading with $\ell$ indices; the corresponding Euler characteristic is the multi-variable Alexander polynomial of the link.

There is also an extension of knot Floer homology to singular links \cite{OSSsingular}.

\subsection{Some basic properties} \label{sec:basic}
We will mostly focus our discussion on $\HFKhat$, and on knots rather than links. (However, the properties below have analogues for the other versions, and for links.)

We start by listing a few symmetries of the knot Floer homology $\HFKhat$:
\begin{itemize}
\item It is insensitive to changing the orientation of the knot;
\item If $m(K)$ denotes the mirror of $K$, we have
$$ \HFKhat_i(K, s) \cong \HFKhat^{-i}(m(K), s),$$
where $\HFKhat^*$ denotes the knot Floer cohomology, related to $\HFKhat_*$ by the universal coefficients formula. Thus, if we use rational coefficients, then the vector space $\HFKhat^i(K, s; \Q)$ is the dual of $ \HFKhat_i(K, s; \Q)$. With $\Z$ coefficients there is an additional Ext term; 
\item We have yet another symmetry:
$$ \HFKhat_i(K, s) \cong \HFKhat_{i-2s}(K, -s).$$
\end{itemize}

There is a K\"unneth formula for the knot Floer homology of connected sums. We state it here for $\Q$ coefficients:
$$ \HFKhat_i(K_1 \# K_2, s; \Q)  \cong \bigoplus_{\substack{i_1+i_2=i\\ s_1+s_2 = s}} \HFKhat_{i_1}(K_1, s_1; \Q) \otimes \HFKhat_{i_2}(K_2, s_2;\Q).$$
Over $\Z$, there is an additional Tor term, as in the usual K\"unneth formula.

Another important property of knot Floer homology is that it admits skein exact sequences. There is an oriented skein exact sequence \cite{Knots}, relating the knot Floer homologies of links that differ at a crossing as follows:
$$
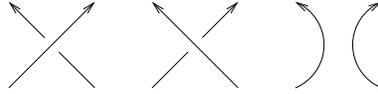
\begin{picture}(0,0)%
\includegraphics{fourlinks.pstex}%
\end{picture}%
\setlength{\unitlength}{2368sp}%
\begingroup\makeatletter\ifx\SetFigFont\undefined%
\gdef\SetFigFont#1#2#3#4#5{%
  \reset@font\fontsize{#1}{#2pt}%
  \fontfamily{#3}\fontseries{#4}\fontshape{#5}%
  \selectfont}%
\fi\endgroup%
\begin{picture}(3920,924)(889,-1573)
\end{picture}%
$$

There is also an unoriented exact sequence \cite{MUnoriented}, relating the knot Floer homologies of a link with those of its two resolutions at a crossing:
$$
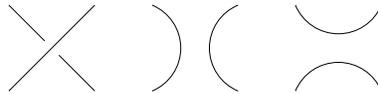
\begin{picture}(0,0)%
\includegraphics{skein.pstex}%
\end{picture}%
\setlength{\unitlength}{2368sp}%
\begingroup\makeatletter\ifx\SetFigFont\undefined%
\gdef\SetFigFont#1#2#3#4#5{%
  \reset@font\fontsize{#1}{#2pt}%
  \fontfamily{#3}\fontseries{#4}\fontshape{#5}%
  \selectfont}%
\fi\endgroup%
\begin{picture}(3920,924)(889,-1573)
\end{picture}%
$$

Yet another exact sequence \cite{CubeResolutions} relates the knot Floer homology of a knot with that of its oriented resolution at a crossing, and with that of its singularization at that crossing:
$$
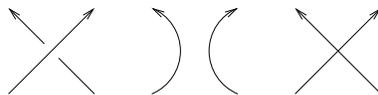
\begin{picture}(0,0)%
\includegraphics{singularskein.pstex}%
\end{picture}%
\setlength{\unitlength}{2368sp}%
\begingroup\makeatletter\ifx\SetFigFont\undefined%
\gdef\SetFigFont#1#2#3#4#5{%
  \reset@font\fontsize{#1}{#2pt}%
  \fontfamily{#3}\fontseries{#4}\fontshape{#5}%
  \selectfont}%
\fi\endgroup%
\begin{picture}(3924,924)(889,-1573)
\end{picture}%
$$

In all three cases, one needs to normalize $\HFKhat$ by tensoring with an additional factor (depending on the number of components of each link). We refer to \cite{Knots, MUnoriented, CubeResolutions} for the details.

\subsection{Calculations}
A knot $K$ is called {\em alternating} if it admits a planar diagram in which the over- and under-passes alternate, as we follow the knot. In the case of alternating knots, knot Floer homology is determined by two classical invariants, the Alexander polynomial $\Delta_K$ and the knot signature $\sigma(K)$. 

\begin{theorem}[Ozsv\'ath-Szab\'o \cite{AltKnots}]
\label{thm:alt}
Let $K \subset S^3$ be an alternating knot with Alexander-Conway polynomial $\Delta_K(q) = \sum_{s \in \Z} a_s q^s$ and signature $\sigma = \sigma(K)$.
Then:
$$ \HFKhat_i(K, s) = \begin{cases}
\Z^{|a_s|} & \text{if} \ i=s + \tfrac{\sigma}{2}, \\
0 & \text{otherwise.}
\end{cases}
$$ 
\end{theorem}
Using the unoriented exact triangle, it can be shown that the same result holds for a more general class of knots, called {\em quasi-alternating} \cite{MOQuasi}. (Most knots with a small number of crossings are quasi-alternating.)

Another class of knots for which $\HFKhat$ is determined by classical invariants (albeit in a different way) is Berge knots, that is, those that can produce a lens space by surgery. (See Section~\ref{sec:surgery} for the definition of surgery.) For example, torus knots are in this class. More generally, one can talk about L-space knots, those that admit a surgery with Heegaard Floer homology ``as simple as possible''; see \cite{OSLens} for the exact definition.

\begin{theorem}[Ozsv\'ath-Szab\'o \cite{OSLens}]
\label{thm:berge}
If $K$ is a Berge knot (or, more generally, an L-space knot), then its Alexander polynomial is of the form
$$\Delta_K(q) = \sum_{j=-k}^k (-1)^{k-j} q^{n_j},$$
for some $k \geq 0$ and integers $n_{-k} < \dots < n_k$ such that $n_{-j} = - n_j$. Furthermore, if we set
$$ \delta_j = \begin{cases}
0 & \text{if} \ j=k \\
\delta_{j+1} - 2(n_{j+1} - n_j) + 1 & \text{if } \ k-j \ \text{is odd} \\
\delta_{j+1}-1 &\text{if } \ k-j > 0 \text{ is even,}
\end{cases}$$
then the knot Floer homology of $K$ has the form
$$ \HFKhat_i(K, s) = \begin{cases}
\Z & \text{if} \ i=n_j \text{ and } s =  \delta_j \text{ for some } j, \\
0 & \text{otherwise.}
\end{cases}
$$ 
\end{theorem}

It is worth mentioning that both for quasi-alternating knots and L-space knots, more is true: the Alexander polynomial and the signature determine not just $\HFKhat$, but the whole full knot complex. See \cite{Petkova, OSLens}.

Hedden \cite{HedCable1, HedWhitehead, HedCable2} and Eftekhary \cite{Eftekhary} studied the knot Floer homology of  cables and Whitehead doubles. Since then, the knot Floer homology of satellites has been further studied using the bordered Floer homology of Lipshitz-Ozsv\'ath-Thurston \cite{LOT}. See also \cite{LevineDouble, Petkova, HomBord} for more recent work in this direction.

Knot Floer homology can also be calculated for many small knots using combinatorial methods. A table with calculations for all non-alternating knots with up to 12 crossings can be found in \cite{BaldwinGillam}. 

\subsection{Geometric applications}
We recall the definition of the genus of a knot:
$$ g(K) = \min \{\text{genus}(F)  \mid F \subset S^3 \mbox{ is an oriented, embedded surface with } \del F = K \}.$$
A well-known property of the Alexander polynomial is that its degree gives a lower bound on the genus \cite{SeifertKnotGenus}. Precisely, if 
$$ \Delta_K(q) = a_0 + a_1(q+q^{-1}) + \dots + a_n(q^{n}+q^{-n}), \ \ a_n \neq 0,$$
then $g(K) \geq n.$

Knot Floer homology strengthens this property, in that it detects the knot genus exactly:
\begin{theorem}[Ozsv\'ath-Szab\'o, Theorem 1.2 in \cite{GenusBounds}]
\label {thm:genusK}
For any knot $K \subset S^3$, we have
$$ g(K) = \max \{s \geq 0 \mid \HFKhat_*(K, s) \neq 0 \}.$$
\end {theorem}

Since the unknot is the unique knot of genus zero, we have
\begin {corollary}
\label{cor:unknot}
If $K \subset S^3$ has the same bigraded knot Floer homology groups $\HFKhat$ as the unknot $U$ (i.e., $\Z$ in bidegree $(0,0)$ and zero otherwise), then $K=U$. 
 \end {corollary}

More generally, the link Floer homology of a link $L \subset S^3$ determines the Thurston norm of the link complement \cite{OSThurstonNorm}.

A knot $K$ is called {\em fibered} if its complement $S^3\setminus K$ fibers over the circle. Another property of the Alexander polynomial is that it provides an obstruction to fiberedness: if $K$ is fibered, then $\Delta_K(q)$ must be monic. Again, knot Floer homology strengthens this property, because it can tell exactly when a knot is fibered:

\begin {theorem}[Ozsv\'ath-Szab\'o, Ghiggini, Ni, Juh\'asz]
\label{thm:fiberedK}
A knot $K \subset S^3$ is fibered if and only if 
$$ \HFKhat_*(K, g(K)) \cong \Z.$$
\end {theorem}

The ``only if" part of the theorem was first proved by Ozsv\'ath-Szab\'o \cite{HolDiskSymp}. Ghiggini proved the ``if'' part for genus one knots \cite{Ghiggini}, and Ni proved it in general \cite{NiFibered}. An alternative proof was given by Juh\'asz \cite{JuhaszDecompose}, using sutured Floer homology.

The only genus one fibered knots are the figure-eight and the two trefoils (right-handed or left-handed). Since their knot Floer homologies can easily be seen to be distinct, we have:
\begin {corollary} [Ghiggini \cite{Ghiggini}]
\label{cor:gh}
Let $E$ be the left-handed trefoil, the right-handed trefoil, or the figure-eight knot. If $K \subset S^3$ has the same bigraded knot Floer homology groups $\HFKhat$ as $E$, then $K=E$. 
\end {corollary}

Knot Floer homology can also be successfully applied to questions of knot concordance. Two knots $K_0$ and $K_1$ are called (smoothly) {\em concordant} if there is a smoothly embedded annulus 
$ A \subset S^3 \times [0,1]$ with $A \cap (S^3 \times \{i\}) = K_i \times \{i\}$ for  $i=0,1.$ A knot concordant to the unknot is called {\em slice}. In fact, there is a notion of {\em slice genus} for a knot:
$$ g_4(K) = \min \{\text{genus}(F)  \mid F \subset B^4 \mbox{ is an oriented, properly embedded surface, } \del F = K \subset S^3 \},$$
and $K$ is slice if and only if $g_4(K)=0$. One reason the slice genus is an interesting quantity is because it gives a lower bound for the {\em unknotting number} $u(K)$ of the knot, that is, the minimum number of crossing changes needed to transform a planar diagram for $K$ into one for the unknot.

One can extract from knot Floer homology an invariant $\tau(K) \in \zz$, which has the property that $\tau(K_1) = \tau(K_2)$ if $K_1, K_2$ are concordant; see \cite{4BallGenus}, \cite{RasmussenThesis}. 
To define $\tau$, one needs more information than the one in $\HFKhat$. One definition (cf. \cite{OST}) can be given in terms of the $\Z[U]$-module $\HFKm$:
$$\tau(K) = -\max \{ s \mid \exists x \in \HFKm_*(K, s), \ U^j x \neq 0 \text{ for all } j \geq 0 \}.$$

The invariant $\tau$ yields an obstruction to two knots being concordant. Further, $\tau$ provides a lower bound on the slice genus of a knot, and hence for the unknotting number:
$$ |\tau(K)| \leq g_4(K) \leq u(K). $$
This allows one to compute the slice genus of various knots. In particular, Ozsv\'ath and Szab\'o used it in \cite{4BallGenus} to give a new proof of a conjecture of Milnor on the slice genus of torus knots (originally proved by Kronheimer and Mrowka using gauge theory \cite{KMMilnor}). 

By using the full knot Floer complex, one can extract additional concordance invariants; see \cite{Hom1, Hom2, HomWu}. Hom \cite{Hom2} applied these ideas to show that the smooth concordance group of topologically slice knots admits a $\Z^{\infty}$ summand.

Using its relation to the three-manifold invariants (outlined in Section~\ref{sec:surgery}), knot Floer homology was successfully applied to questions about surgery. An early example is the constraint on the Alexander polynomial of Berge knots provided by Theorem~\ref{thm:berge}. Another is the work of Ozsv\'ath-Szab\'o \cite{RatSurg}, Wu \cite{Wu} and Ni-Wu \cite{NiWu} on cosmetic surgeries. 

Knot Floer homology has further applications to contact geometry, as it allows the construction of  invariants for Legendrian and transverse knots in $S^3$; see \cite{OST, NOT, LOSS, BVVV}.

\subsection{Distinguishing knots}
A natural question is what knot types can be distinguished by knot Floer homology. From \eqref{eq:alex} we see that if $K_1$ and $K_2$ are distinguished by the Alexander polynomial, then they are also distinguished by knot Floer homology. However, knot Floer homology is a strictly stronger invariant. For example:
\begin {itemize}
\item If $m(K)$ denotes the mirror of $K$, then $\Delta_K = \Delta_{m(K)}$. On the other hand, $\HFKhat(K) \neq \HFKhat(m(K))$ for the trefoil, and for many other knots;
\item If $K_1, K_2$ differ from each other by Conway mutation, then $\Delta_{K_1} = \Delta_{K_2}$. A well-known example of mutant knots, the Conway knot and the Kinoshita-Terasaka knot, have different knot Floer homologies \cite{Mutation}. 
\end {itemize}

As can be seen from the table in \cite{BaldwinGillam}, knot Floer homology is generally an effective invariant for distinguishing between two small knots. Nevertheless, it has its limitations: as mentioned above, in the case of alternating knots, knot Floer homology is determined by the Alexander polynomial and the signature. In particular, we can find examples of different knots with the same knot Floer homology (and, in fact, with the same full knot Floer complex up to filtered homotopy equivalence). The alternating knots $7_4$ and $9_2$ are the simplest such example. 

A related question is what knots $E$ are distinguished from all other knots by knot Floer homology. At present, the only known examples are the four simplest knots: the unknot, the two trefoils, and the figure-eight; cf. Corollaries~\ref{cor:unknot} and \ref{cor:gh}.

\section {The original definition} \label{sec:knot}

We review here the holomorphic curves definition of (several variants of) knot Floer homology, following Ozsv\'ath-Szab\'o \cite{Knots} and Rasmussen \cite{RasmussenThesis}. We will use the more general set-up from \cite{MOS}, allowing multiple basepoints.

\subsection{Heegaard diagrams}

Let $K \subset S^3$ be an oriented knot.

\begin {definition}
\label{def:heegaard}
A (multi-pointed) {\em Heegaard diagram} 
$$(\Sigma, \alphas, \betas, \ws, \zs)$$
for the knot $K$ consists of the following data:
\begin{itemize}
\item A surface $\Sigma \subset S^3$ of genus $g \geq 0$, splitting $S^3$ into two handlebodies $U_0$ and $U_1$, with $\Sigma$ oriented as the boundary of $U_0$;
\item A collection $\alphas=\{\alpha_1, \dots, \alpha_{g+k-1}\}$ consisting of $g+k-1$ pairwise disjoint, simple closed curves on $\Sigma$, such that each $\alpha_i$ bounds a properly embedded disk $D^\alpha_i$ in $U_0$, and the complement of these disks in $U_0$ is a union of $k$ balls $B^\alpha_1, \dots, B^\alpha_k$;
\item A curve collection $\betas = \{\beta_1, \dots, \beta_{g+k-1}\}$ with similar properties, bounding disks $D^\beta_i$ in $U_1$, such that their complement is a union of $k$ balls $B^\beta_1, \dots, B^\beta_k$;
\item Two collections of points on $\Sigma$, denoted $\ws = \{w_1, \dots, w_k\}$ and $\zs = \{z_1, \dots, z_k\}$, all disjoint from each other and from the $\alpha$ and $\beta$ curves.
\end{itemize}
We require that the knot $K$ intersects $\Sigma$ exactly at the $2k$ points $w_i$ and $z_i$, with the intersections being positively oriented at $w_i$ and negatively oriented at $z_i$. Further, we require that the intersection of $K$ with the handlebody $U_0$ consists of $k$ properly embedded intervals, one in each ball $B^\alpha_i$; and similarly that its intersection with $U_1$ consists of $k$ properly embedded intervals, one in each $B^\beta_i$.
\end{definition}

Every knot can be represented by a multi-pointed Heegaard diagram. In fact, if one wishes, it can be represented by a doubly-pointed diagram (i.e., one with $k=1$, so that there is a single $w$ basepoint and a single $z$ basepoint).

\begin{remark}
There is a more general class of Heegaard diagrams for a knot, in which one allows free basepoints; see \cite[Section 4.1]{LinkSurg}. These appear naturally in the context of the link surgery formula discussed in Section~\ref{sec:ls}. 
\end{remark}

\begin {figure}
\begin {center}
\input {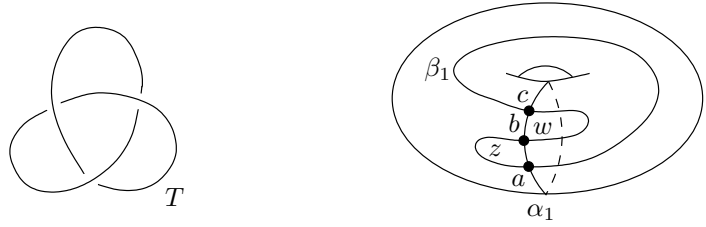}
\caption {The left-handed trefoil knot $T$, and a Heegaard diagram for it.}
\label {fig:trefoil}
\end {center}
\end {figure}

Figure~\ref{fig:trefoil} shows a doubly pointed Heegaard diagram for the trefoil. In general, one can  construct a Heegaard diagram for a knot from a suitable Morse function on the knot complement. There are also more concrete constructions, of which we give a few examples below.

\begin{example}[\cite{AltKnots}] \label{ex:alt}
Suppose we are given a planar diagram for the knot $K$, as the one for the trefoil on the left-hand side of Figure~\ref{fig:alt}. The knot projection $\pi(K)$ is a self-intersecting curve in the plane; it splits the plane into a number of regions $R_0, \dots, R_g$, with $R_0$ being unbounded. Let $\Sigma$ be a boundary of the tubular neighborhood of $\pi(K)$ in $\rr^3$; this is a surface of genus $g$. Draw an alpha curve on $\Sigma$ around each bounded region $R_i$ for $i > 0$. Further, draw a beta curve around each crossing of $K$ as in Figure~\ref{fig:alt}, and an additional beta curve $\beta_g$ as a meridian on $\Sigma$ next to an edge on the boundary of $R_0$. Finally, place two basepoints $w$ and $z$ on each side of $\beta_g$. We get a Heegaard diagram for $K$, as on the right hand side of Figure~\ref{fig:alt}. 

\begin {figure}
\begin {center}
\input {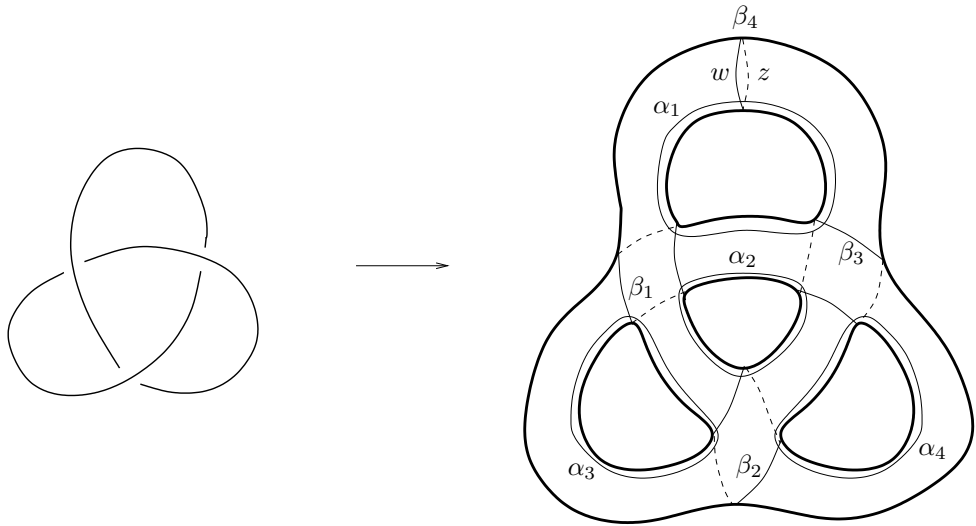}
\caption {A doubly pointed Heegaard diagram associated to a planar projection.}
\label {fig:alt}
\end {center}
\end {figure}

\end{example}

\begin{example} [\cite{RasmussenThesis}]
\label{ex:bridge}
Consider a bridge presentation of the knot, that is, a planar diagram consisting of $2m$ segments $a_1, \dots, a_m, b_1, \dots, b_m$ such that the $a$ curves do not intersect each other, the $b$ curves do not intersect each other, and whenever an $a$ curve crosses a $b$ curve, the $b$ curve is the overpass. See the left hand side of Figure~\ref{fig:bridge} for a bridge presentation of the trefoil, with $m=2$. To a bridge presentation with $m > 1$ bridges we can associate a Heegaard diagram of genus $m-1$, as follows. For each $i=1, \dots, m-1$, draw two shaded disks at the two endpoints of the $a_i$ curve, and identify them; this has the effect of ``adding a handle'' to the plane. Together with the point at infinity, this produces the desired Heegaard surface of genus $m-1$. In the process, each $a_i$ curve (for $i < m$) has become a circle, which we denote by $\alpha_i$. The remaining $a_m$ segment is deleted, and we place the two basepoints $w$ and $z$ at its endpoints. Finally, around $m-1$ of the $b$ curves we draw circles (encircling the handles and/or the basepoints), and denote these by $\beta_1, \dots, \beta_{m-1}$. The result is a Heegaard diagram for $K$, as on the right hand side of Figure~\ref{fig:bridge}. (In the case of the trefoil considered here, this is in fact the same diagram as the one in Figure~\ref{fig:trefoil}.)

\begin {figure}
\begin {center}
\input {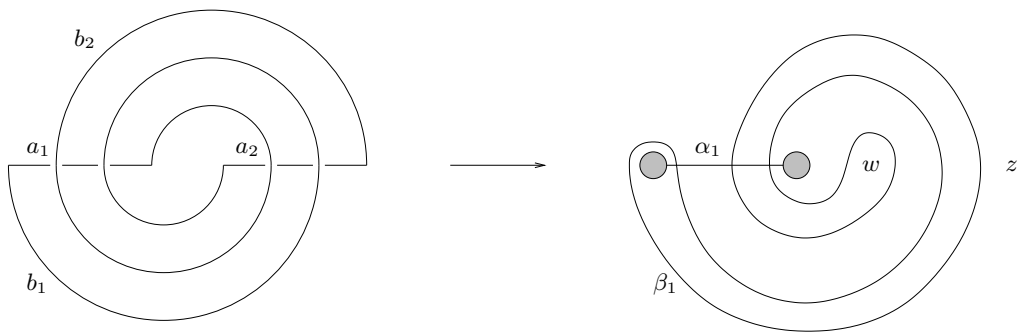}
\caption {A doubly pointed Heegaard diagram associated to a bridge presentation.}
\label {fig:bridge}
\end {center}
\end {figure}

\end{example}

\begin {figure}
\begin {center}
\input {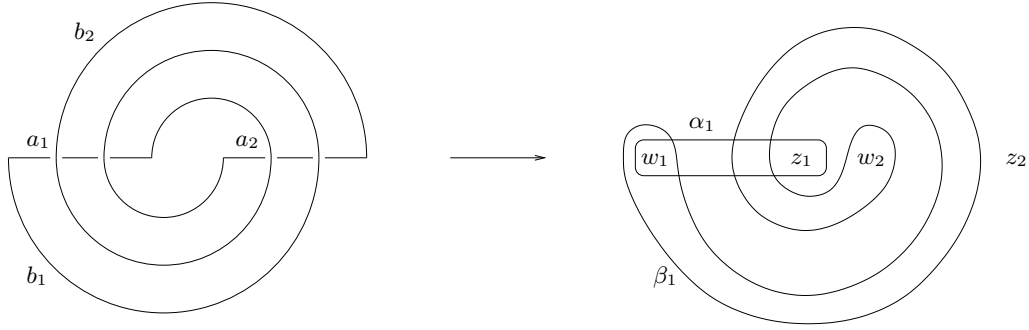}
\caption {A multi-pointed, genus zero Heegaard diagram associated to a bridge presentation.}
\label {fig:planar}
\end {center}
\end {figure}

\begin{example} 
Instead of minimizing the number of basepoints in a Heegaard diagram, we may be interested in minimizing the genus of $\Sigma$. For every knot $K$, we can construct a multi-pointed Heegaard diagram of genus zero as follows. Start with a bridge presentation with $m > 1$, as in the previous example. For each segment $a_i$ for $i < m$, instead of adding a handle, draw a circle $\alpha_i$ around the respective segment, and add basepoints $w_i$ and $z_i$ at the endpoints of $a_i$. Add two extra basepoints $w_m$ and $z_m$ at the endpoints of $a_m$ (but no circle). Lastly, draw beta curves around $b_1, \dots, b_{m-1}$ just as in Example~\ref{ex:bridge}. See Figure~\ref{fig:planar}.
\end{example}

\begin{example}
Another example of Heegaard diagrams are grid diagrams, which will be discussed in detail in Section~\ref{sec:grids}. They are multi-pointed diagrams of genus $1$.
\end{example}

\subsection{Background for the knot Floer complex}
Let $\He = (\Sigma, \alphas, \betas, \ws, \zs)$ be a Heegaard diagram for $K \subset S^3$. Let $g$ be the genus of $\Sigma$ and $k$ the number of basepoints in $\ws$ (or in $\zs$). Let $d = g+k-1$ be the number of alpha curves (which is the same as the number of beta curves).

Starting from this data, we will define the various versions of the knot Floer complex.

Let $\Sigma^{\times d}$ be the Cartesian product of $d$ copies of $\Sigma$. The symmetric group $S_d$ acts on $\Sigma^{\times d}$ by permuting the factors. The quotient is the symmetric product
$$ \Sym^d(\Sigma) := \Sigma^{\times d}/S_d,$$
which is a smooth $2d$-dimensional manifold. Inside of $\Sym^d(\Sigma)$ we consider the half-dimensional submanifolds
$$ \Ta = \alpha_1 \times \cdots \times \alpha_d, \ \ \ \Tb = \beta_1 \times \cdots \times \beta_d$$
obtained by projection from $\Sigma^{\times d}$. We drop the projection from notation for simplicity.

A complex structure on $\Sigma$ induces one on $\Sym^d(\Sigma)$, with respect to which the tori $\Ta, \Tb$ are totally real. (In fact, one can equip $\Sym^d(\Sigma)$ with a symplectic form, such that $\Ta$ and $\Tb$ are Lagrangian; see \cite{Perutz}.) Leaving aside many technicalites, this allows one to define the Lagrangian Floer homology of the pair $(\Ta, \Tb)$.  Roughly, this is the homology of a complex generated by intersection points $\x \in \Ta \cap \Tb$, and whose differential counts pseudo-holomorphic disks in $\Sym^d(\Sigma)$ with boundaries on $\Ta$ and $\Tb$. This kind of construction was first proposed by Floer \cite{FloerLagrangian}, and then developed by various authors \cite{FloerHoferSalamon}, \cite{OhFloer}, \cite{OhTransversality}, \cite{OhFactorize}, \cite{FOOO1}, \cite{FOOO2}.  

In our setting, note that every intersection point $\x \in \Ta \cap \Tb$ consists of an unordered $d$-tuple of points on $\Sigma$, one on each alpha curve and one on each beta curve. We arrange so that the alpha and beta curves intersect transversely. Then, $\Ta \cap \Tb$ is a finite set of points. These will be the generators of the Lagrangian Floer complex, which in this case is called the {\em knot Floer complex}.

If $\x, \y \in \Ta \cap \Tb$ are two intersection points, we denote by $\pi_2(\x, \y)$ the set of relative homotopy classes of disks $u: D^2 \to \Sym^g(\Sigma)$, with $u(-1) = \x, u(1)=\y$, and $u$ taking the lower half of $\del D^2$ to $\Ta$ and the upper half to $\Tb$. Given $\phi \in \pi_2(\x, \y)$, a {\em pseudo-holomorphic representative} for $\phi$ is a map $u$ in the class $\phi$ that satisfies the non-linear Cauchy-Riemann equations with respect to a suitable family of almost complex structures $J=(J_t)_{t \in [0,1]}$ on $\Sym^d(\Sigma)$. (The theory of pseudo-holomorphic curves was initiated by Gromov \cite{Gromov}. See \cite{McDuffSalamon} for an introduction to the subject.) Note that the definition of pseudo-holomorphic depends on $J$; an alternate, synonymous term is {\em J-holomorphic}. However, the space of possible $J$ is contractible, and it turns out that the Floer homology groups will be independent of $J$. 

We denote by $\M(\phi)$ the space of pseudo-holomorphic representatives of $\phi$. Associated to $\phi$ is a quantity $\mu(\phi) \in \zz$, called the {\em Maslov index}. The Maslov index can be calculated using a formula due to Lipshitz \cite{LipshitzCyl}, but this is beyond the scope of this article. (We will explain how $\mu$ can be calculated in more specific examples.) For now, let us note that for generic $J$, the space $\M(\phi)$ is a smooth manifold of dimension $\mu(\phi)$. (In particular, if $\mu(\phi) < 0$, then $\phi$ has no pseudo-holomorphic representatives.) There is an action of $\rr$ on $\M(\phi)$ given by the automorphisms of the domain $D^2$ that fix $1$ and $-1$. Provided that $\phi$ is non-trivial (that is, it is not the class of a constant map), the quotient $\Mh(\phi)= \M(\phi)/\rr$ is smooth of dimension $\mu(\phi)-1$. When $\mu(\phi)=1$, it consists of a discrete set of points. By the general principle of Gromov compactness, $\Mh(\phi)$ is in fact a finite set of points. 

The moduli spaces can be given orientations, depending on some choices; see \cite{HolDisk} for details.\footnote{Orientations have been constructed in \cite{HolDisk, Knots} only for the case of doubly pointed diagrams. It is a ``folklore theorem'' that the same thing can be done for multi-pointed diagrams, but no account of this exists in the literature, except in the particular case of grids \cite{MOST, Gallais}.  For the purposes of this survey, we will assume that orientations can be given, and hence that our complexes are defined over $\zz$. Many sources work over the field $\zz/2$, and then the issue of orientations can be safely ignored. It should be noted that most of the properties of knot Floer homology mentioned in Section~\ref{sec:appl} still hold with coefficients in $\zz/2$.}When $\mu(\phi)=1$, we can then define a signed count of pseudo-holomorphic disks, 
$$\# \Mh(\phi) \in \zz.$$

Each basepoint $v \in \{w_1, \dots, w_k, z_1, \dots, z_k \}$ produces a codimension two submanifold 
$$R_v = \{v\} \times \Sym^{d-1}(\Sigma)$$ inside $\Sym^d(\Sigma)$. By construction, $\Ta$ and $\Tb$ are disjoint from $R_v$. Note also that
$$ \Sym^d(\Sigma \setminus \{v\}) = \Sym^d(\Sigma) \setminus R_v.$$

\begin {figure}
\begin {center}
\input {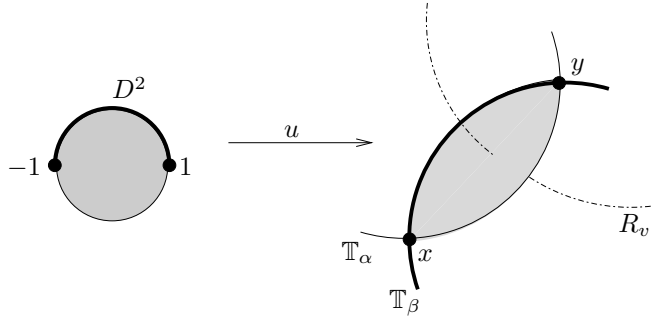}
\caption {A disk $u$ in a class $\phi \in \pi_2(\x, \y)$, and its intersection with $R_v$.}
\label {fig:disk}
\end {center}
\end {figure}

Given intersection points $\x, \y \in \Ta \cap \Tb$, and a class $\phi \in \pi_2(\x, \y)$, we define $n_v(\phi)$ to be the intersection number between $\phi$ and $R_v$. See Figure~\ref{fig:disk}.

\begin {fact}
\label {fact:pos}
(a) (Positivity of intersections) If the class $\phi$ has a pseudo-holomorphic representative, then $n_v(\phi) \geq 0$.

(b) If $n_v(\phi) = 0$, then a generic pseudo-holomorphic representative of $\phi$ does not intersect $R_v$.
\end {fact}

The intersection points $\x \in \Ta \cap \Tb$ admit a bigrading $(M, A)$. The {\em Maslov (homological) grading} $M(\x) \in \zz$ is characterized (up to the addition of a constant) by the formula:
\begin {equation}
\label{eq:Maslov}
 M(\x) - M(\y) = \mu(\phi) - 2 \sum_{i=1}^k n_{w_i} (\phi),
 \end {equation}
where $\phi$ is any relative homotopy class in $\pi_2(\x, \y)$. It can be shown that the right hand side of \eqref{eq:Maslov} is independent of the choice of $\phi$. We should also mention that there is a way of fixing $M$ as an absolute grading in $\zz$, rather than only up to the addition of a constant; this is explained in Section~\ref{sec:minus} below.

The second assignment $A : \Ta \cap \Tb \to \zz$ is called the {\em Alexander grading}, and is uniquely determined by the following two properties:

(i) For any $\x, \y \in \Ta \cap \Tb$ and $\phi \in \pi_2(\x, \y)$, we have
\begin {equation}
\label{eq:Alexander}
 A(\x) - A(\y) = \sum_{i=1}^k n_{z_i}(\phi) - \sum_{i=1}^k n_{w_i}(\phi). 
\end {equation}

(ii) We have
\begin{equation}
\label{eq:alexander}
 \sum_{\x \in \Ta \cap \Tb} (-1)^{M(\x)} q^{A(\x)} = (1- q^{-1})^{k-1} \cdot \Delta_K(q),
 \end{equation}
where $\Delta_K$ is the Alexander-Conway polynomial of the knot $K$. 

\begin {example}
\label{ex:tr}
Consider the Heegaard diagram for the trefoil from Figure~\ref{fig:trefoil}. This has $g=1, k=1, d=1$. While pseudo-holomorphic disks are hard to count in general, when $d=1$ we are in the first symmetric product $\Sym^1 (\Sigma) = \Sigma$, and pseudo-holomorphic disks simply correspond to disks (bigons) on the surface. There are two such disks, one from $a$ to $b$ containing $z$, and one from $c$ to $b$ containing $w$. The corresponding classes $\phi \in \pi_2(a, b), \psi \in \pi_2(c, b)$ have Maslov index $1$. The signed counts $\# \Mh(\phi)$ and $\# \Mh(\psi)$ are $\pm 1$, since the disks are unique in their class. Whether the sign is plus or minus depends on the choices made for  orientations, but different choices will produce the same knot Floer homology. We may assume that 
$$\# \Mh(\phi) = \# \Mh(\psi) =1.$$
Since $d=1$, the divisors $R_z$ and $R_w$ are the points $z$ and $w$ themselves. Thus, we have $n_z(\phi)=n_w(\psi)=1$ and $n_z(\psi) = n_w(\phi) =0$.

The differences in Maslov grading between $a, b, c$ can be computed using  \eqref{eq:Maslov}:
$$M(a) - M(b) = M(b) - M(c) = 1.$$
As we shall see in Section~\ref{sec:minus}, in fact we have: 
 \begin{equation}
 \label{eq:MasTr}
  M(a)=2, \ \ \ M(b) =1, \ \ \ M(c) =0.
  \end{equation}
 
The differences in Alexander grading can be computed using \eqref{eq:Alexander}. The grading is then normalized using the Alexander-Conway polynomial of $T$, which is $\Delta_T(q) = q^{-1}-1+q$. We obtain:
$$ A(a)=1, \ \ \ A(b) =0, \ \ \ A(c) =-1.$$
\end{example}

\subsection{The knot Floer complex, crossing no basepoints} \label{sec:hat} The different flavors of the knot Floer complex have to do with different ways of keeping track of the quantities $n_{z_i}$ and $n_{w_i}$. 
We start by defining the simplest Heegaard Floer complex $\gCFKt(\He)$, in which we only consider disks that do not pass over any basepoints; that is, we restrict attention to classes $\phi \in \pi_2(\x, \y)$ with $n_{w_i}(\phi) = n_{z_i} (\phi)=0$ for all $i$.

As an Abelian group, the complex $\gCFKt(\He)$ is freely generated by intersection points $\x \in \Ta \cap \Tb$. The differential is given by:
\begin {equation}
\label {eq:deltilde}
 \del \x = \sum_{\y \in \Ta \cap \Tb} \sum_{\substack{\phi \in \pi_2(\x, \y) \\ \mu(\phi)=1; \ n_{z_i}(\phi) = n_{w_i}(\phi)=0, \forall i }} \bigl( \# \Mh(\phi) \bigr)  \cdot \y. \end {equation}

By construction, the differential $\del$ decreases the Maslov (homological) grading by one, and keeps the Alexander grading constant. The fact that $\del^2 = 0$ follows from an application of Gromov compactness. The homology of $\gCFKt(\He)$ is denoted $\HFKt(\He)$. The bigrading descends to $\HFKt$, and we denote by $\HFKt_i(\He,s)$ the group in Maslov grading $i$ and Alexander grading $s$.

The situation originally considered in \cite{Knots, RasmussenThesis} was when the Heegaard diagram $\He$ is doubly pointed (that is, $k=1$). In that case $\gCFKt(\He)$ is denoted $\gCFKhat(\He)$, and the homology $\HFKt(\He)$ is denoted $\HFKhat(\He)$. Moreover, we have:

\begin{theorem}[Ozsv\'ath-Szab\'o \cite{Knots}, Rasmussen \cite{RasmussenThesis}]
The isomorphism class of $\HFKhat(\He)$, as a bigraded Abelian group, is an invariant of the knot $K \subset S^3$.
\end{theorem}

Thus, we may denote $\HFKhat(\He)$ by $\HFKhat(K)$. This is the hat version of knot Floer homology, mentioned in Section~\ref{sec:general}. Note that Equation~\eqref{eq:alexander} implies that the Euler characteristic of $\HFKhat$ is  $\Delta_K$, as advertised in Equation~\eqref{eq:alex}.

\begin {example}
In the diagram for the trefoil $T$ discussed in Example~\ref{ex:tr}, both holomorphic disks cross basepoints. Therefore, the differential on $\CFKhat$ is trivial, and we get
$$ \HFKhat_i(T, s) = \begin{cases}
\zz & \text{ for } (s, i) = (-1,0), (0,1), \text{ or } (1,2), \\
0 & \text{ otherwise.}
\end {cases}$$
The trefoil is an alternating knot, and our result is in agreement with Theorem~\ref{thm:alt}.
\end{example}

Let us now go back to the more general set-up, when $\He$ is allowed to have more than two basepoints. This situation was first considered in \cite{MOS}, where it is proved that $\HFKt$ depends on $\He$ only in a mild way. In terms of the number $k$ of $w$ basepoints, the group $\HFKt$ is isomorphic to $2^{k-1}$ copies of $\HFKhat(K)$, with some shifts in degree. Precisely, we have
$$ \HFKt(\He) \cong \HFKhat(K) \otimes V^{\otimes (k-1)},$$
where $V$ is an Abelian group freely generated by an element in bi-degree $(-1,-1)$ and one in bi-degree $(0,0)$. 

\begin{example}
In the genus zero, multi-pointed diagram for the trefoil from Figure~\ref{fig:planar}, the alpha and beta curves intersect each other in six points. This gives six generators for $\gCFKt$, and again there are no disks without basepoints. Thus, $\HFKt$ has rank $6$. Checking the gradings we see that $\HFKt(\He) \cong \HFKhat(T) \otimes V$, as expected.
\end{example}

\begin{remark} 
The reader may wonder why we put the letter $g$ in front of the complexes $\gCFKt$ and $\gCFKhat$, but not in front of their homology. This is because we follow the notation in \cite[Section 11.3]{LOT} and reserve the names $\CFKt$ and $\CFKhat$ for the filtered complexes defined in the next subsection. See Remark~\ref{rem:conv} below for more comments about notation. 
\end{remark}

\subsection{The knot Floer complex, crossing basepoints of one type} \label{sec:minus}The next setting we consider is when we allow the disks to cross the $w$ basepoints, but not the $z$ basepoints. Let us introduce variables $U_i$ to keep track of the basepoints $w_i$, for $i=1, \dots, k$. We define a new  version of the knot Floer complex, $\gCFKm(\He)$, as a module over the ring $\zz[U_1, \dots, U_k]$, freely generated by $\Ta \cap \Tb$, and equipped with the differential
\begin {equation}
\label {eq:delminus}
 \del \x = \sum_{\y \in \Ta \cap \Tb} \sum_{\substack{\phi \in \pi_2(\x, \y) \\ \mu(\phi)=1; \ n_{z_i}(\phi) =0, \forall i }} \bigl( \# \Mh(\phi) \bigr)  \cdot U_1^{n_{w_1}(\phi)} \dots U_k^{n_{w_k}(\phi)} \cdot \y. \end {equation}

This is still bigraded, with each $U_i$ decreasing Maslov grading by $2$ and Alexander grading by $1$. The homology of $\gCFKm(\He)$ is the minus version of knot Floer homology, denoted $\HFKm(K)$. It can be shown that all the variables $U_i$ act the same on homology, and hence $\HFKm(K)$ can be viewed as a $\zz[U]$-module, where $U$ is any of the $U_i$. 
(Of course, we could just start with a doubly pointed Heegaard diagram, and then we would have a single variable $U$ from the beginning.) The isomorphism type of $\HFKm(K)$ is a knot invariant.

There are a couple of related constructions:
\begin{enumerate}[(a)]
\item If instead of using $k$ different variables $U_i$ in the complex, we use a single variable $U$ with exponent $n_{w_1} + \dots + n_{w_k}$ (or, in other words, we set $U_1 = \dots = U_k$), the resulting homology is $\HFKm(K) \otimes V^{\otimes (k-1)}$, where $V$ is the rank two free Abelian group from the previous subsection;
\item If we set a single one of the $U_i$ variables to zero when defining the complex, the resulting homology is $\HFKhat(K)$, with trivial action by the other $U_i$ variables.
\end{enumerate}

Instead of allowing disks to go over $w_i$ and not $z_i$, we could allow them to go over $z_i$ and not $w_i$. We could then use variables $U_i$ to keep track of $n_{z_i}$, and the result would be the same $\HFKm(K)$. However, it is customary to encode this information in another way, using a filtered complex. We let $\CFKt(\He)$ be the complex freely generated (over $\zz$) by $\Ta \cap \Tb$, with differential
\begin{equation}
\label{eq:cfS3}
\del \x = \sum_{\y \in \Ta \cap \Tb} \sum_{\substack{\phi \in \pi_2(\x, \y) \\ \mu(\phi)=1; \ n_{w_i}(\phi) =0, \forall i }} \bigl( \# \Mh(\phi) \bigr)  \cdot \y.
\end{equation}  
Note that the $z_i$ basepoints play no role in this definition, so the knot $K$ disappears from the input. In fact, $\CFKt(\He)$ is a Heegaard Floer complex associated to $S^3$ itself, and its homology (the Heegaard Floer homology of $S^3$) is the homology of a torus, $H_*(T^{k-1})$. We can reintroduce the $z_i$ basepoints by considering the Alexander grading on generators, which defines a filtration on $\CFKt(\He)$, called the {\em knot filtration}. Precisely, the intersection points $\x \in \Ta \cap \Tb$ with $A(\x) \leq j$ form a subcomplex $\F(K, j) \subset \CFKt(\He)$, and we have: 
$$ \dots \subseteq \F(K, j-1) \subseteq \F(K, j) \subseteq \F(K, j+1)  \subseteq \dots$$

The associated graded complex $\bigoplus_j \F(k, j)/\F(k, j-1)$ is $\gCFKt(K)$, and the information in the filtered chain homotopy type of $\CFKt(\He)$ is roughly equivalent to that in $\gCFKm(\He)$. 

One advantage of using the complex $\CFKt(\He)$ is that it helps us fix the absolute Maslov grading. Since we know that the homology of $\CFKt(\He)$ is isomorphic to $H_*(T^{k-1})$ as a relatively graded group, the convention is to set the homological grading so that the top degree element in homology is in degree zero. Together with the relation \eqref{eq:Maslov}, this determines the Maslov grading on generators $\x \in \Ta \cap \Tb$.

The most common situation considered in the literature is for doubly pointed Heegaard diagrams. Then $\CFKt(\He)$ is denoted $\CFKhat(\He)$, and its homology is the hat Heegaard Floer homology of $S^3$, namely $\HFhat(S^3) \cong \zz$ in homological degree zero. The corresponding knot filtration is the one discussed in \cite{Knots}.

\begin{example}
Consider the diagram for the trefoil from Example~\ref{ex:tr}. 
In $\gCFKm(T)$, there is a contribution to the differential from the bigon from $c$ to $b$, which goes over $w$ but not $z$. We get that $\del c = U \cdot b$, so (forgetting the gradings): 
$$\HFKm(T) \cong \zz[[U]] \oplus \zz.$$

If we consider the complex $\CFKhat(\He)$ instead, we have $\del a = b$ and $\del c = \del b = 0$. The homology $\HFhat(S^3)$ is generated by $c$, so we set the Maslov grading of $c$ to be zero. This fixes the Maslov grading of the other generators, and we obtain \eqref{eq:MasTr}.
\end{example}

\subsection{The knot Floer complex, involving all basepoints} \label{sec:both}
Let us now consider the most general situation, in which we allow pseudo-holomorphic disks to cross both types of basepoints. One way of encoding this is to combine the constructions defined by \eqref{eq:delminus} and \eqref{eq:cfS3}. Precisely, we let $\CFKm(\He)$ be the complex freely generated by $\Ta \cap \Tb$ over $\zz[U_1, \dots, U_k]$, equipped with the differential
\begin{equation}
\label{eq:mS3}
\del \x = \sum_{\y \in \Ta \cap \Tb} \sum_{\{\phi \in \pi_2(\x, \y) \mid \mu(\phi)=1 \}} \bigl( \# \Mh(\phi) \bigr) \cdot U_1^{n_{w_1}(\phi)} \dots U_k^{n_{w_k}(\phi)} \cdot \y.
\end{equation}  

Again, the $z_i$ basepoints play no role, and the homology of $\CFKm(\He)$ is a variant of the Heegaard Floer homology of $S^3$, namely $\HF^-(S^3) \cong \zz[U]$ (where each $U_i$ variable acts by $U$). However, the Alexander grading defines a filtration on this complex, which depends on the knot. The filtered chain homotopy type of $\CFKm(\He)$ is a knot invariant. It is usually denoted $\CFKm(K)$, although of course its filtered {\em isomorphism} type depends on the diagram $\He$, not just on the knot.

\begin{remark}
Many homological invariants (for example, the singular homology of a topological space $X$, or the hat knot Floer homology of a knot $K \subset S^3$) are secretly chain homotopy types of free complexes. When the base ring is a PID such as $\zz$, the chain homotopy type of a free complex is determined by its homology; thus, we lose no information by passing to homology. By contrast, when we have a filtered complex, its filtered chain homotopy type is not easily determined by any version of homology. This is why in the case of the more complicated knot invariants discussed here, we can only refer to them as filtered chain homotopy types.
\end{remark}

It is convenient to consider a larger complex than $\CFKm$, denoted $\CFKi$. This is freely generated by $\Ta \cap \Tb$ over the Laurent polynomial ring $\zz[U_1, \dots, U_k, U_1^{-1}, \dots, U_k^{-1}]$, with the differential given by the same formula \eqref{eq:mS3}. The filtered chain homotopy type of $\CFKi$ is again a knot invariant. This is the {\em full knot Floer complex}, previously mentioned in Section~\ref{sec:general}. 

For simplicity, let us restrict to the case of doubly pointed Heegaard diagrams, and write $U$ for $U_1$. In this setting, following \cite{Knots}, we can think of $\CFKi$ as freely generated over $\zz$ by triples 
$$[\x, i, j] , \  \x \in \Ta \cap \Tb, \ i, j \in \zz \text{ with } A(\x) = j-i.$$
The triple $[\x, i, j]$ corresponds to the generator $U^{-i} \x$. Graphically, we represent each generator of $\CFKi$ by a dot in the plane, with $[\x, i, j]$ having coordinates $(i, j)$. Since the action of $U$ decreases Alexander grading by $1$, we can think of the $j$ coordinate as describing the Alexander grading of the generator, whereas the $i$ coordinate describes the (negative of the) $U$ power. The differentials are drawn by arrows. If we have a contribution to the differential from a disk in a class $\phi\in \pi_2(\x, \y)$, note that the change in horizontal coordinate is $-n_w(\phi)$, and the change in vertical coordinate is $-n_z(\phi)$. The Maslov grading is not shown in the picture.

\begin {figure}
\begin {center}
\input {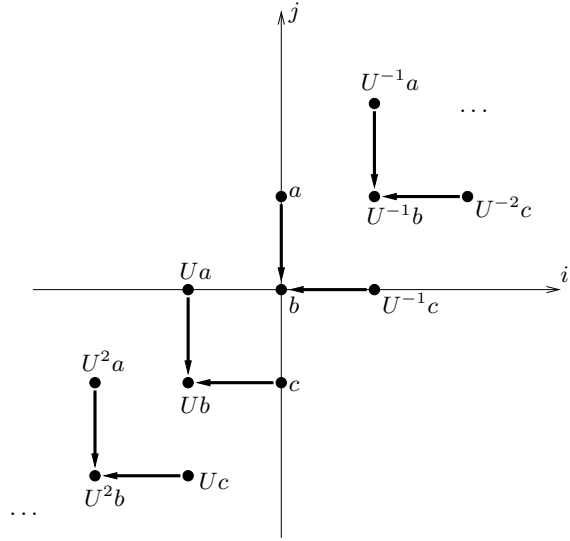}
\caption {The full knot Floer complex for the left-handed trefoil.}
\label {fig:full}
\end {center}
\end {figure}

\begin{example}
For the trefoil $T$ as in Example~\ref{ex:tr}, the full knot Floer complex is drawn in Figure~\ref{fig:full}.
\end{example}

\begin{remark}
Instead of viewing $\CFKi$ as a $\zz$-filtered complex over $\zz[U, U^{-1}]$, we could think of it as a $(\zz \oplus \zz)$-filtered complex over $\zz$, with the two components of the filtration being the coordinates $i$ and $j$.
\end{remark}

From $\CFKi$ we can obtain various auxiliary complexes, by restricting to suitable regions in the $(i, j)$ plane. For example:
\begin{itemize}
\item
 The subcomplex of $\CFKi$ corresponding to triples $[\x, i, j]$ with $i \leq 0$ is exactly $\CFKm$. \footnote{In \cite{Knots}, $\CFKm$ is identified with the subcomplex of $\CFKi$ corresponding to $i < 0$ rather than $i \leq 0$. Of course, the $i < 0$ and $i \leq 0$ complexes are isomorphic (via multiplication by $U$).} 
 \item
 The quotient complex of $\CFKi$ made of triples $[\x, i, j]$ with $i \geq 0$ is denoted $\CFKp.$ 
\item 
The subcomplex of $\CFKp$ made of triples $[\x, i, j]$ with $i=0$ is exactly $\CFhat(S^3)$, with the knot filtration given by the $j$ coordinate;
\item If we consider triples $[\x, i, j]$ with $i=0$ and only allow differentials that preserve the $(i, j)$ grading, we obtain the complex $\CFKhat$.
\end{itemize}

Other auxiliary complexes of this type are the stable complexes that we will discuss in Section~\ref{sec:large}.

\begin{remark}
\label{rem:conv}
There are various notational conventions in the literature with respect to knot Floer complexes. We  followed the notation from \cite[Section 11.3]{LOT}, with $\CFKhat$ and $\CFKm$ denoting filtered complexes whose associated graded are $\gCFKhat$ and $\gCFKm$. In other sources, for example in \cite{LOSS}, the names $\CFKhat$ and $\CFKm$ are used for these associated graded. The original source \cite{Knots} took a mixed approach: the associated graded denoted $\gCFKhat$ here was called $\CFKhat$ in \cite{Knots}, but $\CFKm$ was used to denote the same filtered complex as we did here.   
\end{remark}

\subsection{Links}
The constructions above can be generalized to the case of links \cite{Links}. Let $L \subset S^3$ be a link with $\ell$ components. A Heegaard diagram for $L$ is defined in the same way as for a knot. Note that the minimum number of basepoints is $2\ell$. The generators $\x \in \Ta \cap \Tb$ have a Maslov grading as before. They also have $\ell$ different Alexander gradings, each corresponding to counting $w$ and $z$ basepoints on a link component. 

By analogy with Section~\ref{sec:hat}, we obtain the hat version of link Floer homology, $\HFLhat(L)$, which is an Abelian group with $\ell + 1$ gradings. By analogy with Section~\ref{sec:minus}, we obtain another version, $\HFLm(L)$, in the form of a multi-graded module over the polynomial ring $\zz[U_1, \dots, U_\ell]$. Furthermore, as in Section~\ref{sec:both}, we can consider a more complete invariant, the multi-filtered chain homotopy type of a complex $\CFLm$ over $\zz[U_1, \dots, U_\ell]$.

\section {Grid diagrams}
\label {sec:grids}
The definition of knot Floer homology in the previous section involves counts of pseudo-holomorphic curves. In this section we describe a certain class of Heegaard diagrams, called toroidal grid diagrams, for which the pseudo-holomorphic curve counts become combinatorial. 

The exposition is inspired from the original references \cite{MOS, MOST}.

\subsection {Planar grid diagrams}

Knots (and links) in $S^3$ are usually described in terms of their planar projections. An alternative way to represent them is through grids:

\begin{definition}
A (planar) {\em grid diagram} $G$ of size $n$ is an $n$-by-$n$ grid in the plane, together with $O$ and $X$ markings inside its cells, such that every row and every column contain exactly one $O$ marking and exactly one $X$ marking.
\end{definition}

\begin{definition}
Let $G$ be a grid diagram. Let us trace oriented segments from each $O$ to the $X$ marking in the same row, and from each $X$ marking to the $O$ marking in the same column. We adopt the convention that at each crossing, the vertical segments are overpasses and the horizontal segments are underpasses. We obtain a planar diagram for some oriented link $L$. We then say that the grid diagram $G$ {\em represents} $L$. 
\end{definition}

\begin {figure}
\begin {center}
\begin{picture}(0,0)%
\includegraphics{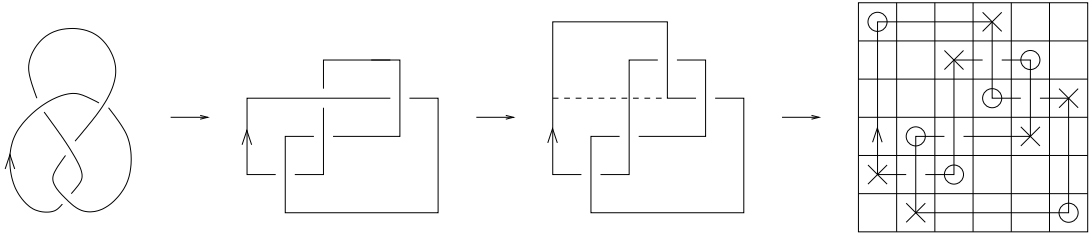}%
\end{picture}%
\setlength{\unitlength}{1579sp}%
\begingroup\makeatletter\ifx\SetFigFont\undefined%
\gdef\SetFigFont#1#2#3#4#5{%
  \reset@font\fontsize{#1}{#2pt}%
  \fontfamily{#3}\fontseries{#4}\fontshape{#5}%
  \selectfont}%
\fi\endgroup%
\begin{picture}(17032,3624)(-6219,-3973)
\end{picture}%

\caption {Transforming a planar projection of the figure-eight knot into a grid diagram.}
\label {fig:constructgrid}
\end {center}
\end {figure}

It is easy to see that every link admits a grid diagram. Indeed, one can start with an ordinary planar projection of $L$, straighten out all segments so that they are either vertical or horizontal, then use small isotopies to place the vertical segments on top of the horizontal ones, and finally draw the grid around the segments. This process is illustrated in Figure~\ref{fig:constructgrid}.

Grid diagrams are equivalent to arc presentations of links, which go back to the work of Brunn \cite{Brunn}; see also \cite{Cromwell, Dynnikov}. The minimal size of a grid needed to present a link $L$ is called the {\em arc index} of $L$. 

In tables, small knots are usually listed in the order of crossing number (the minimal number of crossings in a planar projection). However, an alternate listing can be done in terms of the arc index, as in \cite{JinPark}. It is known that for alternating links, the arc index equals the crossing number plus 2 \cite{BaePark}, while for non-alternating prime links, the arc index is smaller than or equal to the crossing number \cite{JinPark}. 

\subsection{Toroidal grid diagrams} 
Consider a grid diagram representing an oriented link $L$. Let us identify the opposite sides of the square, to obtain a torus. The result is called a {\em toroidal grid diagram}, and can be viewed as a special case of a Heegaard diagram for $L$. Indeed, we let the torus be the Heegaard surface, the $O$ markings be the $w$ basepoints, and the $X$ markings be the $z$ basepoints. Furthermore, the horizontal and vertical lines that form the grid now become circles, and these circles are the alpha and beta curves on the diagram.

\subsection {Combinatorial knot Floer complexes}
Since a toroidal grid diagram $G$ is a particular kind of (multi-pointed) Heegaard diagram, all the different constructions of knot Floer complexes from Sections~\ref{sec:hat}-\ref{sec:both} can be applied to this setting. Note that intersection points $\x \in \Ta \cap \Tb$ correspond to $n$-tuples of points on the grid, with one point on each alpha curve and one point on each beta curve. Thus, the knot Floer complex has exactly $n!$ generators.

What is interesting is that the counts $\# \M(\phi)$ of pseudo-holomorphic curves become very concrete.  Indeed, it is proved in \cite{MOS} that generic, index one pseudo-holomorphic curves are in one-to-one correspondence with {\em empty rectangles} on the grid. 

\begin {figure}
\begin {center}
\begin{picture}(0,0)%
\includegraphics{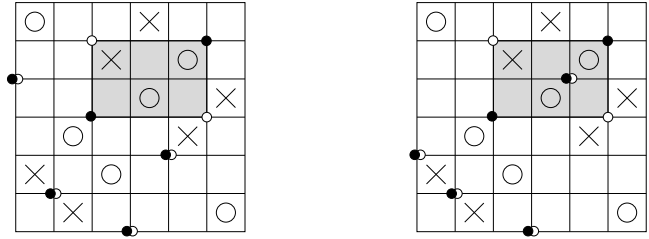}%
\end{picture}%
\setlength{\unitlength}{1579sp}%
\begingroup\makeatletter\ifx\SetFigFont\undefined%
\gdef\SetFigFont#1#2#3#4#5{%
  \reset@font\fontsize{#1}{#2pt}%
  \fontfamily{#3}\fontseries{#4}\fontshape{#5}%
  \selectfont}%
\fi\endgroup%
\begin{picture}(10050,3687)(1063,-4036)
\end{picture}%

\caption {We show an empty rectangle (contributing to the differential in $\CFKi$) shaded on the left, and a non-empty rectangle shaded on the right. In each picture, the generator $\x$ is shown as a $5$-tuple of black dots, and the generator $\y$ as a $5$-tuple of white dots. (Note that the top edge is identified with the bottom edge, so there are actually components of $\x$ and $\y$ there too, which we did not draw. Similarly, there are components of $\x$ and $\y$ on the right edge.)}
\label {fig:rectangles}
\end {center}
\end {figure}

\begin{definition}
Let $G$ be a toroidal grid (viewed as a Heegaard diagram), and let $\x, \y \in \Ta \cap \Tb$. We view $\x$ and $\y$ as $n$-tuples of points on the grid. A {\em rectangle} from $\x$ to $\y$ is an embedded rectangle in $G$, such that its bottom edge is an arc on an alpha curve from a point of $\x$ (on the left) to a  point of $\y$ (on the right), and its top edge is an arc on an alpha curve from a point of $\y$ (on the left) to a point of $\x$ (on the right); furthermore, we assume that the remaining $n-2$ components of $\x$ coincide with the remaining $n-2$ components of $\y$.

The rectangle is called {\em empty} if it does not contain any of the remaining $n-2$ components of $\x$ (or $\y$) in its interior. See Figure~\ref{fig:rectangles}.  
\end{definition}

Note that rectangles live on the torus, not on the plane, so they can wrap around the edges of the planar grid diagram.

Let $\EmptyRect(\x, \y)$ be the set of empty rectangles from $\x$ to $\y$. Each rectangle $r \in \EmptyRect(\x, \y)$ has an associated relative homotopy class $\phi \in \pi_2(\x, \y)$, and the quantities $n_{w_i}(\phi)$ and $n_{z_i}(\phi)$ are either $0$ or $1$, according to whether the corresponding basepoint is or is not inside the rectangle. Given that on the grid the basepoints are marked by $X$ and $O$, it is customary to write $\Xs_i(r)$ for $n_{w_i}(\phi)$ and $\Os_i(r)$ for $n_{z_i}(\phi)$.

Moreover, to each rectangle $r \in \EmptyRect(\x, \y)$ one can associate a sign $\eps(r) \in \{\pm 1\}$, which is meant to represent the orientation of the respective pseudo-holomorphic curve. We refer to  \cite{MOST, Gallais} for the exact formula for $\eps$.

With this in mind, the knot Floer complexes from Sections~\ref{sec:hat}-\ref{sec:both} become purely combinatorial. For example, $\CFKt(G)$ is a free Abelian group generated by the $n!$ possible $n$-tuples of points $\x$, with the differential:
$$\del \x = \sum_{\y \in \Ta \cap \Tb} \sum_{\substack{r \in \EmptyRect(\x, \y) \\ \Os_i(r) = \Xs_i(r) = 0,\ \forall i}} \eps(r) \cdot \y.$$

One can give a completely combinatorial proof of the invariance of knot Floer complexes, using grids; see \cite{MOST}.

The grid diagram method was implemented on the computer \cite{BaldwinGillam, Droz} and used to calculate the knot Floer homology of knots with small arc index.  

Nevertheless, it should be noted that making the differentials in the knot Floer complex combinatorial comes at the price of greatly increasing the number of generators (compared with other Heegaard diagrams). For example, in Section~\ref{sec:knot} we computed the different knot Floer complexes of the trefoil $T$ using the Heegaard diagram from Figure~\ref{fig:trefoil}, with $3$ generators. By contrast, the arc index of the trefoil is $5$, so the smallest grid diagram representing $T$ has size $5$ and thus $5! = 120$ generators. 

In general, grid diagrams are useful for computing the knot Floer homology of small knots. If one is interested in special (infinite) classes of knots, other methods may be more helpful. For example, for alternating knots, the Heegaard diagrams constructed in Example~\ref{ex:alt} yield a hat knot Floer complex with no non-trivial differentials, and hence with the minimum possible number of generators. These diagrams were used in \cite{AltKnots} to prove Theorem~\ref{thm:alt}.

\section{The cube of resolutions}
\label{sec:cube}

We now turn to a different combinatorial formulation of knot Floer homology, developed by Ozsv\'ath and Szab\'o in \cite{CubeResolutions}. This is based on constructing a cube of resolutions involving singular links.

\subsection{Definition}
\begin {figure}
\begin {center}
\input {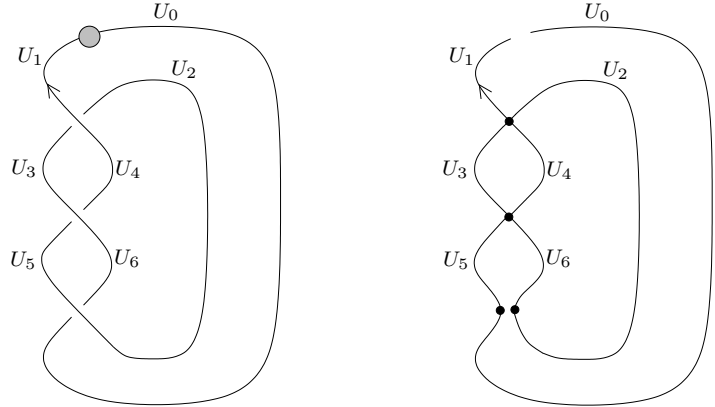}
\caption {A braid diagram for the left-handed trefoil (with the top leftmost edge split in two at the gray dot), and a complete resolution of this diagram (after one crossing was smoothed and two were singularized). We show the variables $U_i$ corresponding to each edge.
}
\label {fig:braid}
\end {center}
\end {figure}

Let $K \subset S^3$ be a knot. We start with a braid presentation $\K$ of $K$, as on the left hand side of Figure~\ref{fig:braid}. Let $n$ be the number of crossings in $\K$. We cut the top leftmost edge in two at a distinguished point (marked by a gray dot), so that there are now $2n+1$ of edges (arcs between crossings, or between a crossing and the gray dot) in the diagram. Let $E = \{e_0, e_1, \dots, e_{2n}\}$ be the set of these edges, such that $e_0$ is the edge starting at the distinguished point, according to the orientation of the knot.  

We can resolve each crossing $p$ in $\K$ in two ways: either by taking the oriented resolution at $p$, and marking a point on each of the two resulting arcs (this is called a {\em smoothing}), or by replacing the crossing with a valence four intersection point between the two arcs (this is called a {\em singularization}). If we do one of these two opertaions at each crossing, the result is called a {\em complete resolution} of $\K$. An example is shown on the right hand side of Figure~\ref{fig:braid}.

Let $\Ring$ be the polynomial ring $\zz[U_0, \dots, U_n]$, where each variable corresponds to an edge. We consider also the ring $\Ring[t]$, where we adjoin an extra variable $t$. To each complete resolution $S$ of $\K$ we associate two ideals $L_S, N_S \subset \Ring[t]$ as follows. 

The ideal $L_S$ is generated by elements $L(p)$, one for each crossing $p$ that is singularized in $S$. If $a(p)$ and $b(p)$ are the outgoing edges from $p$, and $c(p)$ and $d(p)$ are the incoming edges, we set
\begin{equation}
\label{eq:L}
 L(p) = t \cdot (U_{a(p)} + U_{b(p)}) - (U_{c(p)} + U_{d(p)}).
 \end{equation}

The ideal $N_S$ is generated by several elements $N(W)$, one for each subset $W$ of vertices in the graph associated to $S$. Given such a subset $W$, let $|W|$ be the number of smoothed vertices plus twice the number of singular vertices in $W$. Let $\Out(W)$ be the set of outgoing edges from vertices $W$, and $\In(W)$ be the set of incoming edges to vertices in $W$. If $W^c$ denotes the complement of $W$, we define the element 
\begin{equation}
\label{eq:N}
 N(W)= t^{|W|} \cdot \prod_{e_i \in \Out(W) \cap \In(W^c)} U_i - \prod_{e_i \in \In(W) \cap \Out(W^c)} U_i.
 \end{equation}
 
\begin{example}
\label{ex:reso}
Consider the complete resolution $S$ on the right hand side of Figure~\ref{fig:braid}. The ideal $L_S$ is
$$ L_S = (tU_1 + tU_2 - U_3 - U_4, tU_3 + tU_4 - U_5 - U_6).$$
The ideal $N_S$ has $2^4=16$ generators, corresponding to all possible subsets $W$. For example, when $W$ consists of the two singular points, we have $N(W) = t^4 U_1U_2 - U_5U_6.$ If we add to $W$ the rightmost of the two vertices at the smoothing, we get  $t^5U_1 -  U_5.$ Some of these $16$ elements are generated by the others, so in the end we obtain:
$$ N_S = (tU_5-U_0, tU_6 - U_2, t^5U_1 - U_5, t^2 U_1U_2 - U_3 U_4, t^2 U_3U_4 - U_5 U_6).$$ 
\end{example}

Going back to the general situation, we define an algebra $\A(S)$ as the quotient 
$$ \A(S) := \Ring[t]/(L_S + N_S).$$

It is not hard to see that $\A(S)$ is zero when $S$ is disconnected.

\begin {figure}
\begin {center}
\input {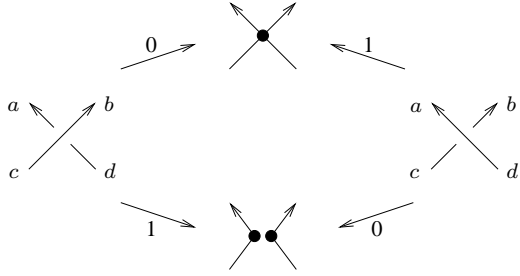}
\caption {Resolutions at a crossing}
\label {fig:res}
\end {center}
\end {figure}

We now construct a cube of resolutions. Let $c(\K)$ be the set of crossings in $\K$. If $p \in c(\K)$ is a positive crossing, we define its $0$-resolution to be its singularization at $p$, and its $1$-resolution to be its smoothing. If $p$ is a negative crossing, we do the opposite. (See Figure~\ref{fig:res}.) Given a map $I: c(\K) \to \{0, 1\}$, we can form the complete resolution $S_I(\K)$ by resolving all the crossings according to $I$. Let

$$ C(\K) = \bigoplus_{I: c(\K)\to \{0,1\}} \A(S_I(\K)).$$

We turn this into a chain complex by equipping it with a differential $D$ as follows. Suppose $I, J : c(\K) \to \{0,1\}$ differ at a single crossing $p \in c(\K)$, where we have $I(p)=0$ and $J(p)=1$. We define a map
$$ D_{I < J} : \A(S_I(\K)) \to \A(S_J(\K))$$
by the following procedure. If $p$ is a positive crossing in $\K$, so that $S_I(\K)$ is singular at $p$ and $S_J(\K)$ is smoothed at $p$, then one can observe that $\A(S_J(K))$ is a quotient of $\A(S_I((\K))$. Then $D_{I < J}$ is called an {\em unzip} homomorphism, and is the natural quotient map. If $p$ is a negative crossing, then $D_{I < J}$ is called a {\em zip} homomorphism, and is induced by multiplication by $t \cdot U_a - U_d$, where $a, b, c, d$ are the edges near $p$ drawn as in Figure~\ref{fig:res}.

The differential
$$ D: C(\K) \to C(\K)$$
is the sum over all maps $D_{I < J}$ as above. It can be shown that $D^2=0$.

Let $\zz[t^{-1}, t]]$ be the ring of semi-infinite Laurent power series in $t$; that is, an element of $\zz[t^{-1}, t]]$ is a formal series $\sum_{n \in \zz} a_n t^n$, with $a_n = 0$ for $n \ll 0$. By tensoring the complex $C(\K)$ with $\zz[t^{-1}, t]]$ we obtain a complex $C(\K)[t^{-1}, t]]$.

\begin{theorem}[Ozsv\'ath-Szab\'o \cite{CubeResolutions}] \label{thm:cube}
If $\K$ is a braid presentation of a knot $K$ as above, then we have isomorphisms:
\begin{align*}
\HFKhat(K) \otimes_{\zz} \zz[t^{-1}, t]] &\cong H_*( C(\K)[t^{-1}, t]]/U_0),\\
\HFKm(K) \otimes_{\zz} \zz[t^{-1}, t]] &\cong H_*( C(\K)[t^{-1}, t]]).
\end{align*}
\end{theorem}

Theorem~\ref{thm:cube} gives a combinatorial construction of $\HFKhat$ and $\HFKm$, the versions of knot Floer homology that involve going over only one kind of basepoints; compare Sections~\ref{sec:hat} and \ref{sec:minus}. One can also recover the Alexander and Maslov gradings from the cube of resolutions. On the other hand, it is not known whether the full knot Floer complex from Section~\ref{sec:both} (that involves going over all basepoints) can be obtained from this construction. Further, there is no complete (combinatorial) proof of invariance of knot Floer homology that uses the description above; see however \cite{GilmoreCube} for some progress in this direction.

Let us mention the main ideas that go in the proof of Theorem~\ref{thm:cube}. Unlike in the case of grids, one does not simply count pseudo-holomorphic disks starting from a single Heegaard diagram for $K$. Rather, one uses the third skein exact sequence mentioned in Section~\ref{sec:basic}, and applies it to all possible resolutions and crossings. By a procedure developed in \cite{BrDCov}, iterating the exact sequence yields a spectral sequence, whose $E_1$ page is the cube of resolutions complex $C(\K)$ and that converges to knot Floer homology. Finally, the use of twisted coefficients (the variable $t$)    guarantees that the $E_1$ page lives in a single diagonal line, with respect to the two gradings. This implies that there is no room for higher differentials, so the spectral sequence collapses after taking the $D=d_1$ differential.

\subsection{Comparison with Khovanov-Rozansky homology}

The cube of resolutions above is reminiscent of the Khovanov-Rozansky triply graded homology \cite{KR2}, a knot homology theory whose Euler characteristic is the HOMFLY-PT polynomial. We shall not 
present the definition of the Khovanov-Rozansky homology $H_{\kr}(K)$ in detail, but let us give a sketch. One starts with a braid presentation $\K$ and considers all complete resolutions $S_I(\K)$, just as above. (In \cite{KR2}, one draws singularizations as trivalent graphs with a wide edge, but this is just a cosmetic difference.) To each complete resolution $S_I(\K)$ one associates a complex $C_{\kr}(S_I(\K))$, and then defines
$$ C_{\kr}(\K) = \bigoplus_{I: c(\K) \to \{0,1\}} H_*\bigl (C_{\kr}(S_I(\K))\bigr),$$
with a differential $D_{\kr}$ composed of zip and unzip homomorphisms from $C_{\kr}(S_I(\K))$ to $C_{\kr}(S_J(\K))$, where $I < J$ differ at a single crossing $p$, just as before. Then $H_{\kr}(K)$ is the homology of the complex $(C_{\kr}(\K), D_{\kr}).$

The definition of $C_{\kr}(\K)$ differs from that of $C(\K)$ in the following two ways:
\begin{enumerate}[(a)]
\item In $C(\K)$ we work over the ring $\Ring[t]$, whereas the complex $C_{\kr}(\K)$ is defined over the simpler ring $\Ring = \zz[U_0, \dots, U_{2n}]$, without the $t$ variable;
\item If $S$ is a complete resolution, the homology $H_*(C_{\kr}(S))$ is not a quotient algebra like $\A(S)$, but can live in several gradings. Precisely, $C_{\kr}(S)$ is a Koszul complex, the tensor product of mapping cone complexes of the form 
$$\Ring \xrightarrow{L(p)} \Ring \ \ \text{ and } \ \ \Ring \xrightarrow{N(\{p\})} \Ring,$$ 
where $L(p)$ and $N(\{p\})$ are as in \eqref{eq:L} and \eqref{eq:N}, setting $t=1$. Note that, unlike in $\A(S)$, here we only consider elements $N(W)$ when $W$ consists of a single point $\{p\}$. Thus, the elements $L(p)$ are linear, and $N(\{p\})$ are either linear or quadratic (depending on whether the vertex $p$ is two-valent or four-valent). 
\end{enumerate}

As mentioned in the introduction, Dunfield, Gukov and Rasmussen \cite{DGR} conjectured the existence of a spectral sequence with $E_2$ page the Khovanov-Rozansky homology $H_{\kr}(K)$, and converging to knot Floer homology. A proposal for where this spectral sequence might come from is given in \cite{Spectral}. One sets $t=1$ in the constructions from \cite{CubeResolutions}. Iterating the skein exact sequence yields a spectral sequence converging to knot Floer homology, but now the sequence does not collapse after the $E_1$ page. Still, its $E_1$ page is a complex formed by a cube of resolutions, where at each complete resolution $S$, instead of $\A(S)$ we have a graded $\Ring$-module $\T(S)$. This is easier to describe in the case when $S$ is connected; then $\T(S)$ is the direct sum 
\begin{equation}
\label{eq:torN}
 \bigoplus_i \Tor^{\Ring}_i(\Ring/L_S, \Ring/N_S),
 \end{equation}
where $L_S$ and $N_S$ are defined just as in the previous subsection, except that we have set $t=1$ so that they are now ideals of $\Ring$.

By comparison, in the Khovanov-Rozansky cube, if $S$ is connected then the homology $H_{\kr}(S_I(\K))$ can be shown to be isomorphic to
\begin{equation}
\label{eq:torQ}
 \bigoplus_i \Tor^{\Ring}_i(\Ring/L_S, \Ring/Q_S),
 \end{equation}
where $Q_S$ differs from $N_S$ in that it is only generated by the elements $N(W)$ with $W=\{p\}$ consisting of a single vertex.

It is conjectured in \cite{Spectral} that the $E_1$ page of the $t=1$ spectral sequence above is isomorphic to the Khovanov-Rozansky complex. If true, this would imply the Dunfield-Gukov-Rasmussen conjecture. 

The main stumbling block is an algebraic problem, that of finding an isomorphism between the Tor groups \eqref{eq:torN} and \eqref{eq:torQ}. Some evidence for the existence of such an isomorphism is presented in \cite{Spectral}. 

\begin{example}
Consider the complete resolution $S$ on the right hand side of Figure~\ref{fig:braid}, studied in Example~\ref{ex:reso}. When we set $t=1$, we have:
\begin{align*}
 L &= (U_1 + U_2 - U_3 - U_4, U_3 + U_4 - U_5 - U_6), \\
 Q &= (U_5-U_0, U_6 - U_2,  U_1U_2 - U_3 U_4, U_3U_4 - U_5 U_6), \\
 N &= Q + (U_1 - U_5).
 \end{align*}
We leave it as an exercise to check that the Tor groups \eqref{eq:torN} and \eqref{eq:torQ} are isomorphic in this example.
\end{example}

\begin{remark}
For any $S$, note that $Q_S$ is a subset of the ideal $N_S$. However, it is shown in \cite{Spectral} that  the desired isomorphism between \eqref{eq:torN} and \eqref{eq:torQ} cannot always come from the natural quotient map. A more useful fact may be that $N_S$ is an ideal quotient of $Q_S$, as proved by Gilmore \cite{GilmoreIdeal}.
\end{remark}

\section {Surgery formulas}
\label{sec:surgery}

In this section we outline the relation between knot Floer homology and the Heegaard Floer invariants of three- and four-manifolds. The relation is expressed in terms of surgery formulas, which connect the knot Floer complexes to the invariants of surgeries on the knot.

\subsection{Surgery} \label{sec:surg}
Let us briefly review how three- and four-manifolds can be expressed in terms of links in $S^3$. Good introductory textbooks about this subject are \cite{SavelievBook} and \cite{GompfStipsicz}.

Let $p, q \in \zz$ be two relatively prime numbers, and let $K \subset S^3$ be a knot. The result of $p/q$ {\em Dehn surgery} along $K$ is the three-manifold
$$ S^3_{p/q}(K) = (S^3 - \nu(K)) \cup_{\phi} (S^1 \times D^2).$$
Here $\nu(K)$ is a tubular neighborhood of the knot, and $S^1 \times D^2$ is attached to its boundary by a diffeomorphism $\phi: S^1 \times \del D^2 \to \del \nu(K)$ taking the meridian $* \times \del D^2$ to a curve in the homology class $p[\mu] + q[\lambda]$, where $\lambda$ and $\mu$ are the longitude and meridian of the knot, respectively. The ratio $p/q \in \Q \cup \{\infty\}$ is called the surgery coefficient.

When $q=1$, the surgery is called {\em integral surgery} (or {\em Morse surgery}). In that case, there is an induced four-dimensional cobordism $W$ from $S^3$ to $S^3_p(K)$, given by attaching a two-handle:
$$ W = \bigl( S^3 \times [0,1] \bigr)  \cup_{\psi} (D^2 \times D^2),$$
where $\psi$ takes $\del D^2 \times D^2$ to $\nu(K)$. 

Not every $3$-manifold can be expressed as Dehn surgery on a knot, but a theorem of Lickorish and Wallace \cite{Lickorish, Wallace} says that every (closed, oriented) $3$-manifold $Y$ can be expressed by integral surgery along a link $L \subset S^3$; i.e., by doing integral surgery on each link component. The collection of integral surgery slopes describes a framing $\Lambda$ of the link, and there is an induced cobordism from $S^3$ to $Y$, given by several two-handle attachments.

Lastly, consider a closed, oriented four-manifold $X$. Any such $X$ can be expressed as the union of a zero-handle, several one-handles, two-handles, and three-handles, and a four-handle. This decomposition can be represented graphically by a Kirby diagram; see \cite{GompfStipsicz} for details. Most of the intricacy in four-manifold topology comes from the two-handles, so for simplicity we will focus our discussion on cobordisms consisting of two-handle attachments (i.e., those coming from integral surgeries on links).

\subsection {Heegaard Floer theory} 
\label {sec:hf}
We now give a quick outline of the construction of the Heegaard Floer invariants. For more details, we refer the reader to the expository papers \cite{SurveyOS}, \cite{OSBudapest1}, \cite{OSBudapest2}, \cite{McDuffSurvey}, and to the original articles (referenced below). 

Heegaard Floer theory originated in the papers \cite{HolDisk}, \cite{HolDiskTwo} of Ozsv\'ath and Szab\'o, in which they defined new invariants of 3-manifolds in the form of homology theories $\HFhatold, \HFminus, \HFplus, \HFinf$. These are the various flavors of {\em Heegaard Floer homology}. They are modules over the polynomial ring $\zz[U]$. 

The definition of Heegaard Floer homology is very similar to that of knot Floer homology. It starts with a Heegaard diagram representing a closed, oriented three-manifold $Y$. A Heegaard diagram for $Y$ is a collection $\He=(\Sigma, \alpha, \beta, \ws)$ with the same properties as in Definition~\ref{def:heegaard}, except that $\Sigma$ is a surface inside of $Y$, and we do not have $z$ basepoints. For simplicity, we will only consider the case when $k=1$, so that there is a single basepoint $w$ and the genus $g$ of $\Sigma$ equals the number of alpha curves and the number of beta curves. (However, Heegaard Floer homology can also be defined when there are multiple $w$ basepoints; see \cite{Links}.)

As in Section~\ref{sec:knot}, we consider the symmetric product $\Sym^g(\Sigma)$, and take the Lagrangian Floer homology of the tori $\Ta=\alpha_1 \times \dots \times \alpha_g$ and $\Tb=\beta_1 \times \dots \times \beta_g$. There are several versions. The simplest is the {\em hat} Heegaard Floer homology $\HFhat(Y)$. This is the homology of a complex $\CFhat(\He)$ generated by intersection points $\x \in \Ta \cap \Tb$, with differential
$$ \del \x = \sum_{\y \in \Ta \cap \Tb} \sum_{\substack{\phi \in \pi_2(\x, \y) \\ \mu(\phi)=1; \ n_{w}(\phi) =0}} \bigl( \# \Mh(\phi) \bigr)  \cdot \y.$$
The analogue of the decomposition of $\HFKhat$ into Alexander gradings is a decomposition of $\HFhat(Y)$ according to $\spinc$ structures:
$$ \HFhat(Y) = \bigoplus_{\ss \in \spinc(Y)} \HFhat\circ(Y, \ss).$$
Here, a $\spinc$ structure is a lift of the frame bundle of $Y$ to a principal $\spinc(3)$-bundle; more concretely, $\spinc$ structures are in (non-canonical) one-to-one correspondence with elements of the second cohomology $H^2(Y; \zz)$.

The {\em minus} complex $\CFminus(Y)$ is generated by $\Ta \cap \Tb$ over $\zz[U]$, with differential\footnote{This definition applies to manifolds $Y$ with $b_1(Y)=0$. When $b_1(Y) > 0$, there is an admissibility condition for Heegaard diagrams, which cannot be satisfied for all $\spinc$ structures at a time. Rather, the $\CFminus$ complex is defined in each $\spinc$ structure starting from a diagram admissible for that structure.}
$$ \del \x = \sum_{\y \in \Ta \cap \Tb} \sum_{\substack{\phi \in \pi_2(\x, \y) \\ \mu(\phi)=1}} \bigl( \# \Mh(\phi) \bigr)  \cdot U^{n_w(\phi)} \y.$$

The {\em infinity} and {\em plus} complexes are obtained from $\CFminus$ by
$$ \CFinfty = U^{-1} \CFminus := \CFminus \otimes_{\zz[U]} \zz[U, U^{-1}]$$and $$ \CFplus = \CFinfty/\CFminus.$$

Let $\circ \in \{ \widehat{\; \;}, -, \infty, +\}$ denote any of the four flavors. The homology of $\CF^\circ(\He)$ is denoted $\HF^\circ(Y)$, and is a three-manifold invariant. We have a decomposition of $\HF^\circ(Y)$ according to $\spinc$ structures. There is also a homological grading, which may take values in $\zz$ or in a cyclic group $\zz/d$ (where $d$ depends on the $\spinc$ structure). For simplicity, we will ignore the grading in our exposition.

\begin{example}
When $Y=S^3$, there is a unique $\spinc$ structure, and we have:
$$ \HFhat(S^3) = \zz, \ \ \HFminus(S^3) = \zz[U], \ \ \HFinfty(S^3) = \zz[U, U^{-1}],$$
$$ \HFplus(S^3) = \zz[U^{-1}] := \zz[U, U^{-1}]/\zz[U].$$
\end{example}

Next, suppose we have a four-dimensional cobordism $W$ between three-manifolds $Y_0$ and $Y_1$. Let $W$ have a (four-dimensional) $\spinc$ structure $\tt$ that restricts to $\ss_0$ on $Y_0$ and $\ss_1$ on $Y_1$. In \cite{HolDiskFour}, Ozsv\'ath and Szab\'o construct a cobordism map:
$$ F^\circ_{W, \tt} : \HF^\circ(Y_0, \ss_0) \to \HF^\circ(Y_1, \ss_1).$$

The construction involves decomposing $W$ into handles. The maps associated to one- and three-handles are given by concrete formulas. The maps associated to two-handles are more complicated, being defined in terms of counts of pseudo-holomorphic triangles in $\Sym^g(\Sigma)$ (with boundaries on three Lagrangian tori). 

If $X$ is a closed four-manifold with $b_2^+(X) \geq 2$, by decomposing $X$ into two suitable cobordisms and combining the map on $\HFminus$ on one cobordism with the map on $\HFplus$ on the other cobordism, one can define a {\em mixed} Heegaard Floer invariant for $X$. This is conjecturally the same as the Seiberg-Witten invariant of $X$ \cite{Witten}, which has important applications in four-dimensional topology.

\subsection {Large surgeries}
\label{sec:large}
The simplest surgery formula involves integral surgeries on a knot $K \subset S^3$ where the slope $p \in \zz$ is very large. In this case $S^3_p(K)$ is a rational homology sphere with $H_1 \cong H^2$ of order $p$. It admits $p$ different $\spinc$ structures, which can be identified with the elements of $\zz/p$ in a natural way; see for example \cite[Section 2.4]{IntSurg}.

Let $\He$ be a doubly pointed Heegaard diagram for $K$. Recall that in Section~\ref{sec:both} we defined a full knot Floer complex $\CFKinfty(K)$, represented by a diagram in the $(i, j)$ plane as in Figure~\ref{fig:full}. Recall also that by restricting to various subsets of the plane, we can define auxiliary complexes. 

\begin {figure}
\begin {center}
\input {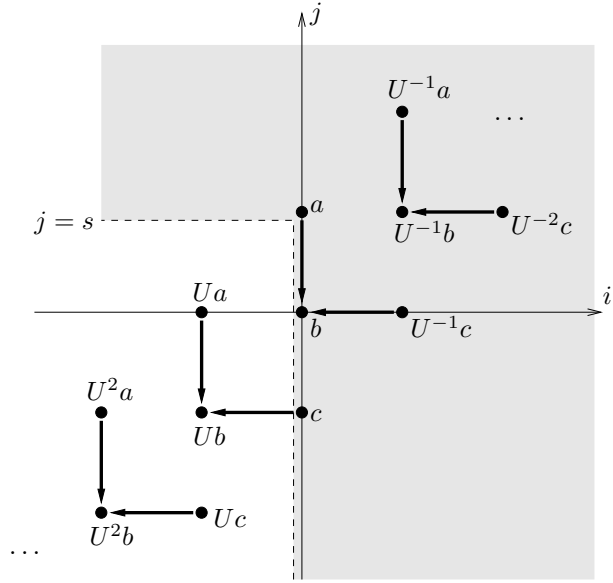}
\caption {The stable knot Floer complex $A_s^+$ for the trefoil is generated by the elements from the shaded region. We draw here the case $s=1$.}
\label {fig:stable}
\end {center}
\end {figure}

For each $s \in \zz$, we define the (plus) {\em stable knot complex} $A_s^+$ to be the quotient complex of $\CFKinfty(K)$ generated by triples $[\x, i, j]$ with $\x \in \Ta \cap \Tb, A(\x) = j-i$ and
$$ \max(i, j-s) \geq 0.$$ 
See Figure~\ref{fig:stable} for an example. The corresponding subcomplex corresponding to the region $\max(i, j-s) < 0$ is denoted $A_s^-$, and the subquotient corresponding to $\max(i, j-s)=0$ is denoted $\hat{A}_s$. Finally, let $A_s^{\infty} = \CFKinfty$, forgetting the filtration.

\begin{theorem}[Ozsv\'ath-Szab\'o \cite{Knots}, Rasmussen \cite{RasmussenThesis}]
\label{thm:large}
Let $\circ \in \{ \widehat{\; \;}, -, \infty, +\}$ be a flavor of Heegaard Floer homology. For $p \gg 0$ and any $s \in \zz$ with $|s| \leq p/2$,  there is an isomorphism
$$ \HF^\circ(S^3_p(K), [s]) \cong H_*(A_s^\circ).$$
Here, $[s] \in \zz/p$ is the corresponding $\spinc$ structure on $S^3_p(K)$. 
\end{theorem}

In other words, the Heegaard Floer homologies of $S^3_p(K)$ ``stabilize'' as $p \to \infty$, and the answer is given by the homology of the corresponding stable knot Floer complex. The proof of Theorem~\ref{thm:large} involves a direct comparison between the stable knot complexes $A_s^\circ$ and the Heegaard Floer complexes $\CF^\circ$ associated to a suitable Heegaard diagram of the surgered manifold.
 
 \begin{example}
 \label{ex:stableT}
Consider the left-handed trefoil with the full knot Floer complex from Figure~\ref{fig:full}. As can be seen in Figure~\ref{fig:stable}, for $s \geq 1$ the complexes $A_s^+$ are all the same, being generated by the elements $U^i \x$, with $\x \in \{a, b, c\}$ and $i \leq 0$. Their homology is generated by the classes $[U^{i+1}a-U^{i}c]$ for $i \leq 0$. We get 
$$\HFplus(S^3_p(T), [s]) \cong H_*(A_s^+) \cong \zz[U^{-1}] \text{ for } s> 0.$$  
When $s=0$, we have an additional generator $Ua$, and the homology of $A_0^+$ is
$$\HFplus(S^3_p(T), [0]) \cong H_*(A_0^+) \cong \zz[U^{-1}] \oplus \zz,$$
with $\zz[U^{-1}]$ coming from $Ua-c$ and $\zz$ from $c$.
When $s \leq -1$, the complexes $A_s^+$ are all isomorphic to each other, via multiplication by powers of $U$, and the homology is once again
$$\HFplus(S^3_p(T), [s]) \cong H_*(A_s^+) \cong \zz[U^{-1}] \text{ for } s< 0.$$  
\end{example}

\begin{remark}
 Theorem~\ref{thm:large} admits a generalization to large surgeries on links, \cite[Theorem 10.1]{LinkSurg}.
\end{remark}

\subsection {The knot surgery formula} Next, consider the case of an arbitrary integral surgery on $K$. 
For $p \in \zz$ nonzero, the manifold $S^3_p(K)$ is a rational homology sphere, again having $|p|$ $\spinc$ structures, in one-to-one correspondence with the elements of $\zz/|p|$. When $p = 0$, the manifold $S^3_0(K)$ has $b_1=1$ and infinitely many $\spinc$ structures, in one-to-one correspondence with the elements of $\zz$.

Consider the stable complex $A_s^+$ from the previous subsection, associated to the region $\max(i, j-s) \geq 0$. It admits two quotient complexes associated to the regions $j \geq s$ and $i \geq 0$. Looking back at Section~\ref{sec:both}, we see that in the complex $B^+$ associated to $i \geq 0$ we ignore the $z$ basepoint; in fact, this is exactly the complex $\CFplus$ for $S^3$. Consider the natural quotient map
$$v_s^+ : A_s^+ \to B^+.$$

Similarly, in the complex $\widetilde B^+$ associated to $j \geq s$ we ignore the $w$ basepoint. Thus, $\widetilde B^+$ is also a $\CFplus$ complex for $S^3$, but corresponding to a Heegaard diagram where we use the $z$ basepoint instead of $w$. Nevertheless, we can find a chain homotopy equivalence between $\widetilde B^+$ and $B^+$. By composing the quotient map $A_s^+ \to \widetilde B^+$ with this homotopy equivalence, we obtain a map
$$h_s^+: A_s^+ \to B^+ .$$

Let us denote $B_s^+ = B^+$ for all $s$. We construct a {\em surgery complex} $\XX^+(p)$ as the mapping cone complex of the map
$$ \bigoplus_{s \in \zz} A_s^+ \longrightarrow \bigoplus_{s\in \zz} B_s^+, \ \ \ (s, \x) \mapsto (s, v_s^+(\x)) + (s+p, h_s^+(\x)).$$
Here, the first component $s$ in a pair $(s, \x)$ indicates the index of the respective direct summand, e.g., $(s+p, h_s^+(\x))$ means that we look at $h_s^+(\x) \in B_{s+p}^+$. 

For example, the complex $\XX^+(2)$ is of the form
\[ \xymatrixcolsep{.7pc}
\xymatrix{
\dots \ar[rrd] & A_{-2}^+\ar[rrd] \ar[d] & A_{-1}^+ \ar[rrd] \ar[d] & A_0^+ \ar[rrd] \ar[d] & A_1^+ \ar[rrd] \ar[d] & A_2^+  \ar[d] & \dots \\
\dots & B_{-2}^+ & B_{-1}^+ & B_0^+ & B_1^+ & B_2^+ & \cdots 
} \]
with the vertical maps being $v_s^+$ and the diagonal maps being $h_s^+$. 

Note that $\XX^+(p)$ naturally splits as a direct sum of two complexes $\XX_i^+(p), \ i \in \zz/p$, where $\XX_i^+(p)$ consists of $A_s^+$ and $B_s^+$ with $s \equiv i \pmod p.$

We can now state the knot surgery formula. Its proof involves using skein exact sequences to relate the Heegaard Floer homology of $S^3_p(K)$ to that of large surgeries on $K$, and then applying Theorem~\ref{thm:large} to make the connection with knot Floer complexes.

\begin{theorem}[Ozsv\'ath-Szab\'o \cite{IntSurg}]
\label{thm:IntSurg}
Let $p$ be an integer, and $K \subset S^3$ an oriented knot. Let $W_p(K)$ be the cobordism (two-handle attachment) from $S^3$ to $S^3_p(K)$ induced by surgery on $K$. 

\noindent $(a)$ There is an identification of the $\spinc$ structures on $S^3_p(K)$ with $\zz/|p|$ such that for any $i \in \zz/|p|$, we have an isomorphism
$$ \HFplus(S^3_p(K), i) \cong H_*(\XX_i^+(p)).$$

\noindent $(b)$ For $s\in \zz$, there is an identification of the $\spinc$ structures on $W_p(K)$ such that the map induced on homology by the inclusion $B^+_s \hookrightarrow \XX_i^+(p)$ produces the cobordism map
$$ F_{W, s}^+: \HFplus(S^3) \to \HFplus(S^3_p(K), [s]).$$
\end{theorem}

\begin{remark}
Even though the complex $\XX^+(p)$ involves two infinite direct sums, these can be truncated so that we are left with finite direct sums. Here, by ``truncating'' we mean deleting an acyclic subcomplex or quotient complex; this process yields a new complex with the same homology. The procedure is illustrated in the two examples below. 
\end{remark}

\begin{example}
Consider again the left-handed trefoil, as in Example~\ref{ex:stableT}. The result of $+1$ surgery on $T$ is the Brieskorn sphere $-\Sigma(2,3,7)$. Let us compute $\HFplus$ of this manifold using Theorem~\ref{thm:IntSurg}. The complexes $A_s^+$ and $B^+$ are described in Example~\ref{ex:stableT}. Furthermore, from Figure~\ref{fig:stable} we see that the map $v_s^+ : A_s^+ \to B^+$ is a homotopy equivalence for $s \geq 1$. Similarly, the map $h_s^+: A_s^+ \to B^+$ is the identity for $s \leq -1$. The complex $\XX^+(1)$ has the form
$$ \input{plus1.pstex_t}$$
where the thicker arrows indicate homotopy equivalences. The shaded area to the right (composed of the complexes $A_s^+$ and $B_s^+$ for $s > 0$) forms a subcomplex of $\XX^+(1)$, denoted $\XX^+(1)\{>0\}$. This subcomplex admits a filtration (indicated by the dashed vertical lines) such that the associated graded decomposes into mapping cone complexes $A_s^+ \to B_s^+$, with the respective map being a homotopy equivalence. These complexes are acyclic, and hence the whole subcomplex $\XX^+(1)\{ > 0\}$ is also acyclic. (We are using here a well-known principle in homological algebra, that under some mild conditions, if we have a filtered complex with trivial homology for its associated graded, then the complex itself is also acyclic.) Similarly, the shaded area to the right forms a subcomplex $\XX^+(1)\{ < 0\}$ which is acyclic because the filtration indicated by diagonal dashed lines produces an acyclic associated graded. Quotienting out the two shaded complexes, we conclude that $\XX^+(1)$ has the same homology as $A_0^+$; the latter was calculated in Example~\ref{ex:stableT}, and we obtain:
$$ \HFplus(S^3_1(T)) \cong \zz[U^{-1}] \oplus \zz.$$
\end{example}

\begin{example}
Let us now consider $-1$ surgery on the left-handed trefoil, which yields the Poincar\'e sphere $\Sigma(2,3,5)$. The complex $\XX^+(-1)$ has the form
$$ \input{minus1.pstex_t}$$
In this case, the shaded regions on the left and right represent two acyclic quotient complexes. We deduce that $\XX^+(-1)$ has the same homology as the subcomplex in the middle, composed of $A_0^+, B_{-1}^+$ and $B_1^+$. The complex $A_0^+$ and the maps $h_0^+, v_0^+$ can be read  off Figure~\ref{fig:stable}. A short calculation shows that
$$ \HFplus(S^3_{-1}(T)) \cong \zz[U^{-1}].$$ 
\end{example}

Observe that Theorem~\ref{thm:IntSurg} is phrased in terms of the plus version of Heegaard Floer homology. The same result holds for the hat version, with the complexes $A_s^+$ replaced by $\hat{A}_s$. There is an analogous result for the minus version, but in that case (for technical reasons, explained in \cite[Section 8.1]{LinkSurg}) one needs to make the following modifications: 
\begin{itemize}
\item Instead of $\HFminus$ one needs to consider the completion of this $\zz[U]$-modules with respect to the $U$ variable. This completion is denoted $\HFMinus$, and yields only a minor loss of information. For rational homology spheres, $\HFminus$ can be recovered from the completed theory;
\item Similarly, one needs to complete $A_s^-$ with respect to $U$;
\item Instead of the direct sums in the construction of $\XX^+(p)$, one needs to take direct products to construct a similar complex $\XX^-(p)$.
\end{itemize}

In a different direction, Theorem~\ref{thm:IntSurg} was generalized to Dehn surgeries on knots with rational surgery slope \cite{RatSurg}. The Heegaard Floer homology $\HFplus(S^3_{p/q})$ is identified with the homology of a mapping cone complex similar to $\XX^+(p)$, except that we use $q$ copies of each $A_s^+$ in the top row.

\subsection {The link surgery formula and four-manifolds} \label{sec:ls}
Another generalization of the knot surgery formula is given in \cite{LinkSurg}. It applies to integral surgeries on arbitrary links $L \subset S^3$. This result is phrased in \cite{LinkSurg} in terms of the completed minus version $\HFMinus$, using direct products to form infinite complexes. (Also, it uses $\zz/2$ rather than $\zz$ coefficients.) Instead of a single mapping cone complex we need to consider a whole hypercube of complexes, where at the vertices we have the stable Floer complexes associated to $L$ and all its sublinks. This hypercube is called the {\em link surgery complex.} Along the edges of the hypercube we have maps similar to $v_s^+$ and $h_s^+$, but one needs to be particularly careful about the choices of chain homotopy equivalences that relate Heegaard Floer complexes for the same geometric object. (Recall that, for knots, we needed such an equivalence between $B^+$ and $\widetilde{B}^+$.) Further, along the diagonals of the hypercube we need to introduce chain homotopies between compositions of the edge maps, then chain homotopies between the new maps, and so on. We refer to \cite{LinkSurg} for the exact formulation.

Part (b) of Theorem~\ref{thm:IntSurg} can also be generalized to links. There is a natural cobordism from surgery on a sublink $L' \subset L$ to surgery on the whole link $L$, given by two-handle attachments along the complement $L-L'$. The cobordism map on $\HFMinus$ induced by this cobordism is exactly the map on homology induced by the inclusion of a subcomplex into the surgery complex; see \cite[Section 11]{LinkSurg} for details.

As mentioned in Section~\ref{sec:surg}, every three-manifold can be obtained by surgery along a link in $S^3$.  Thus, its Heegaard Floer homology can be expressed in terms of Floer complexes associated to links in $S^3$. In turn, these complexes can be described combinatorially using grid diagrams, as in Section~\ref{sec:grids}. For the surgery complex one needs a bit more than complexes, namely maps (and chain homotopies) relating these complexes. With more work, one can show that  the maps and homotopies can also be described combinatorially \cite{MOT}. The result is that the Heegaard Floer homology groups of three-manifolds are algorithmically computable (if we use $\zz/2$ coefficients). 

It was also mentioned in Section~\ref{sec:surg} that four-manifolds can be represented by Kirby diagrams, and that the main intricacy there comes from two-handle attachments. The link surgery formula from \cite{LinkSurg} allows one to describe the maps induced by two-handle attachments in terms of the surgery complex. By using grid diagrams, this description can be made combinatorial. Building on these ideas, one can show that the mixed invariants of closed four-manifolds (with $\zz/2$ coefficients) are algorithmically computable \cite{MOT}. (Unfortunately, the algorithms are very complicated and not yet suitable for practical computations.) See \cite{ECM} for a survey of these developments.

\subsection{Speculations} A topic of current interest is extending the Khovanov and Khovanov-Rozansky homologies \cite{Khovanov, KR1, KR2} to invariants of $3$-manifolds and perhaps $4$-manifolds. Witten \cite{FiveBranes} has made a gauge-theoretic proposal in this direction. There is also some work by Khovanov and Qi \cite{KhQi} aimed at categorifying the Witten-Reshetikhin-Turaev invariants of $3$-manifolds.  

A natural question is whether the surgery formulas for knot Floer homology described in this section can serve as a model for similar constructions in Khovanov-Rozansky homology. Let us list a few ingredients that appear in the surgery formulas from \cite{IntSurg, LinkSurg}:
\begin{enumerate}[(a)]
\item One uses the full knot Floer complex to get the stable complexes $A_s^+$, and there are two different kinds of maps $v_s^+$ and $h_s^+$ relating $A_s^+$ to the Heegaard Floer complex for $S^3$.
\item There is a symmetry between the maps $v_s^+$ and $h_s^+$. For example, the codomain of the former is a complex that uses the $w$ basepoint, and the codomain of the latter is a complex that uses  the $z$ basepoint.
\item When working with links, one uses the full link Floer complex over the module $\zz[U_1, \dots, U_{\ell}]$. In particular, recall that the hat version of link Floer homology categorifies the multi-variable Alexander polynomial.
\end{enumerate}

We discuss each of these points in turn.

With respect to (a), if we look at the hat version, the analogs of $v_s^+$ and $h_s^+$ produce two spectral sequences from $\HFKhat(K)$ to $\HFhat(S^3) \cong \zz$. Thus, the minimum that one needs for a knot surgery formula is two such spectral sequences. In the $\sl(2)$ Khovanov homology there is only one such spectral sequence, due to Lee \cite{Lee}. However, in the triply-graded Khovanov-Rozansky homology from \cite{KR2}, there are two spectral sequences. They were constructed by Rasmussen in \cite{RasDiff}, and correspond to the $d_1$ and $d_{-1}$ differentials conjectured in \cite{DGR}. Therefore, the triply graded homology looks like a natural candidate for a theory with surgery formulas.

The analog of (b) in Khovanov-Rozansky homology would be a symmetry that interchanges the $d_1$ and $d_{-1}$ differentials. This symmetry was conjectured in \cite{DGR}, but its existence has not yet been proved.

For (c), one would need to develop a variant of Khovanov-Rozansky homology for links that categorifies some kind of multi-variable HOMFLY-PT polynomial (with one variable for each link component). This is rather mysterious, since no such natural polynomial is known.

Finally, let us mention the following issue. The result of $+13$ surgery on the torus knot $T(7,2)$ is the lens space $L(13,4)$. The knot $T(7,2)$ is alternating and has ``thin'' triply-graded homology, so by analogy with the Heegaard Floer case, we expect the rank of the presumed HOMFLY-PT homology of $L(13,4)$ to be $13$. If this homology behaves well with respect to orientation reversal, the answer for $-L(13, 4)$ would also have rank $13$. On the other hand, $-L(13,4)$ is $+13$ surgery on $T(3,4)$. If an analog of Theorem~\ref{thm:berge} held then we would get that the HOMFLY-PT polynomial of $T(3,4)$ (with the $a$-grading collapsed) has only $0$ and $\pm 1$ coefficients, but this is not the case. This shows that if a HOMFLY-PT homology for three-manifolds exists, then its properties must be somewhat different from those of Heegaard Floer homology.

\bibliographystyle{amsalpha} 
\bibliography{biblio} 

\end{document}